\documentclass[12pt,reqno,twoside]{amsart}
\usepackage{amsmath,amsthm,amsfonts,amssymb,setspace,verbatim,cite}

\usepackage[english]{babel}
\usepackage[utf8]{inputenc}
\usepackage{times}
\usepackage[T1]{fontenc}
\usepackage{mathrsfs}
\usepackage{mathabx,thmtools}
\usepackage{framed}

\DeclareMathAlphabet{\mathpzc}{OT1}{pzc}{m}{it}

\usepackage{enumitem}

\usepackage[top=15mm,bottom=25mm,left=20mm,right=25mm,includeheadfoot,headsep=8mm,footnotesep=5mm]{geometry}

\usepackage[hypertexnames=false]{hyperref} 
\hypersetup{
pdfencoding=auto,
pageanchor=true,
plainpages=true,
pdftoolbar=true,
pdfmenubar=true,
pdfpagemode=UseOutlines,
pdffitwindow=false,
pdfstartview={FitH},
pdfnewwindow=true,
colorlinks=true,
linkcolor=black,
citecolor=black,
filecolor=black,
urlcolor=black,
linktoc=all,
breaklinks=true,
}

\usepackage{enumitem}

\usepackage{accents}
\DeclareMathAccent{\wwtilde}{\mathord}{largesymbols}{"65}
\DeclareMathSymbol{\widetildesym}{\mathord}{largesymbols}{"65}
\newcommand\lowerwidetildesym{%
  \text{\smash{\raisebox{-1.3ex}{%
    $\widetildesym$}}}}
\newcommand\wtilde[1]{%
  \mathchoice
    {\accentset{\displaystyle\lowerwidetildesym}{#1}}
    {\accentset{\textstyle\lowerwidetildesym}{#1}}
    {\accentset{\scriptstyle\lowerwidetildesym}{#1}}
    {\accentset{\scriptscriptstyle\lowerwidetildesym}{#1}}
}

\usepackage{accents}
\DeclareMathSymbol{\widehatsym}{\mathord}{largesymbols}{"62}
\newcommand\lowerwidehatsym{%
  \text{\smash{\raisebox{-1.3ex}{%
    $\widehatsym$}}}}
\newcommand\what[1]{%
  \mathchoice
    {\accentset{\displaystyle\lowerwidehatsym}{#1}}
    {\accentset{\textstyle\lowerwidehatsym}{#1}}
    {\accentset{\scriptstyle\lowerwidehatsym}{#1}}
    {\accentset{\scriptscriptstyle\lowerwidehatsym}{#1}}
}

\declaretheorem[name=Theorem,numberwithin=section]{theorem}
\declaretheorem[name=Lemma,sibling=theorem]{lemma}
\declaretheorem[name=Corollary,sibling=theorem]{corollary}

\declaretheorem[style=definition,name=Hypothesis,sibling=theorem]{hypothesis}
\declaretheorem[style=definition,name=Remark,sibling=theorem]{remark}

\declaretheorem[style=definition,name=Example,qed=$\blacktriangle$,sibling=theorem]{example}

\numberwithin{equation}{section}

\makeatletter
\def\th@plain{%
  \thm@notefont{\bfseries}
  \itshape 
}

\def\th@definition{%
  \normalfont 
  \thm@notefont{\bfseries}
}
\makeatother

\renewcommand{\thefootnote}{\fnsymbol{footnote}}

\newcommand{\FootCorAuth}
{\footnote{Corresponding author.}}

\newcommand{\makepapertitle}{
\pdfbookmark[1]{Title page}{Title_page}
\thispagestyle{empty}\vspace*{8mm}
\begin{spacing}{1.1}
\LARGE\MakeUppercase{\papertitle}
\end{spacing}\vspace*{10mm}
}

\newcommand{\paperfirstauthor}{%
{\large\sc\firstauthor\FootCorAuth}\\[2mm] 
{\rm\firstaddress \\ 
\firstemail}\\[8mm]
}

\newcommand{\papersecondauthor}{%
{\large\sc\secondauthor}\\[2mm] 
{\rm\secondaddress \\ 
\secondemail}\\[8mm]
}

\newcommand{\makeabstract}{%
\begin{minipage}{0.9\textwidth}
{\small {\sc Abstract.}
 \paperabstract
}
\end{minipage}\vfill
}

\newcommand{\MakeFirstPageTwoAuthors}{
\begin{center}
  \makepapertitle
  \paperfirstauthor
  \papersecondauthor
  \makeabstract
  \begin{minipage}[t]{0.3\textwidth}
   \raggedleft {\bf Date (revised final version):} 
  \end{minipage}\hspace{0.01\textwidth}
  \begin{minipage}[t]{0.63\textwidth}
    \today\\ (submitted on March 09, 2015; accepted on March 08, 2016)
  \end{minipage}\\[4mm]
\noindent%
  \begin{minipage}[t]{0.3\textwidth}
    \raggedleft {\bf Running head:} 
  \end{minipage}\hspace{0.01\textwidth}
  \begin{minipage}[t]{0.63\textwidth}
    \runninghead
  \end{minipage}\\[4mm]
\noindent%
  \begin{minipage}[t]{0.3\textwidth}
    \raggedleft {\bf How to cite:}  
  \end{minipage}\hspace{0.01\textwidth}
  \begin{minipage}[t]{0.63\textwidth} 
    J. Math. Anal. Appl. {\bf 440} (2016), no.~1, 323--350.\\[1mm]
     \footnotesize{\url{http://dx.doi.org/10.1016/j.jmaa.2016.03.028}}
  \end{minipage}\\[4mm]
\noindent%
  \begin{minipage}[t]{0.3\textwidth}
    \raggedleft {\bf License:}  
  \end{minipage}\hspace{0.01\textwidth}
  \begin{minipage}[t]{0.63\textwidth} 
    \copyright{} \the\year. This manuscript version is made available under the 
     \href{http://creativecommons.org/licenses/by-nc-nd/4.0/}{CC-BY-NC-ND 4.0 license}
  \end{minipage}\\
\bigskip
\vfill
\end{center} 

\renewcommand{\thefootnote}{}
\footnotetext[2]{2010 {\it Mathematics Subject Classification:} \thesubjclass}
\footnotetext[3]{{\it Key words and phrases:} \thekeywords}
\setcounter{footnote}{0}
\renewcommand{\thefootnote}{{\bf\,\alph{footnote}\alph{footnote}\alph{footnote}\,}}

\markboth{{\sc\firstauthor\ and \secondauthor}}{{\sc\runninghead}}

\setcounter{page}{0}
\newpage
}

\newcommand{\papertitle}%
{Characterization of self-adjoint extensions for discrete symplectic systems}

\newcommand{\runninghead}%
{Self-adjoint extensions for discrete symplectic systems}

\newcommand{\firstauthor}%
{Petr Zem{\'{a}}nek}

\newcommand{\firstaddress}
{Department of Mathematics and Statistics, Faculty of Science, Masaryk University \\
Kotl{\'{a}}{\v{r}}sk{\'{a}} 2, CZ-61137 Brno, Czech Republic}

\newcommand{\firstemail}%
{E-mail: zemanekp@math.muni.cz}

\newcommand{\secondauthor}%
{Stephen  Clark}

\newcommand{\secondaddress}
{Department of Mathematics \& Statistics, 101 Rolla Building,\\
Missouri University of Science and Technology, Rolla, MO 65409-0020, USA}

\newcommand{\secondemail}%
{E-mail: sclark@mst.edu}

\newcommand{\paperabstract}%
{All self-adjoint extensions of minimal linear relation associated with the discrete symplectic system are characterized. 
Especially, for the scalar case on a finite discrete interval some equivalent forms and the uniqueness of the given expression 
are discussed and the Krein--von Neumann extension is described explicitly. In addition, a limit point criterion for symplectic 
systems is established. The result partially generalizes even a classical limit point criterion for the second order 
Sturm--Liouville difference equations.}

\newcommand{\thekeywords}%
{Discrete symplectic system; linear relation; self-adjoint extension; Krein--von Neumann extension; uniqueness; limit point 
criterion.}

\newcommand{\thesubjclass}%
{{\it Primary\/} 47A06; {\it Secondary\/} 47A20; 39A70; 47B39; 39A12.}

\newcommand{\submittedto}%
{Journal of Mathematical Analysis and Applications}

\newcommand{\e}{\mathrm{e}}

\newcommand{\Ac}{\mathcal{A}}
\newcommand{\Bc}{\mathcal{B}}
\newcommand{\Cc}{\mathcal{C}}
\newcommand{\Dc}{\mathcal{D}}
\newcommand{\Fc}{\mathcal{F}}
\newcommand{\Gc}{\mathcal{G}}

\newcommand{\Ic}{\mathcal{I}}
\newcommand{\tIc}{\wtilde{\Ic}}
\newcommand{\Jc}{\mathcal{J}}

\renewcommand{\Mc}{\mathcal{M}}

\newcommand{\Rc}{\mathcal{R}}
\newcommand{\Sc}{\mathcal{S}}
\newcommand{\tSc}{\wtilde{\Sc}}

\newcommand{\Vc}{\mathcal{V}}
\newcommand{\tVc}{\wtilde{\Vc}}

\newcommand{\Wc}{\mathcal{W}}

\DeclareMathAlphabet{\msfsl}{U}{eus}{m}{n}
\newcommand{\Ff}{\msfsl{F}}
\newcommand{\Gf}{\msfsl{G}}
\newcommand{\Mf}{\msfsl{M}}
\newcommand{\Lf}{\msfsl{L}}
\newcommand{\Vf}{\msfsl{V}}

\newcommand{\mL}{\mathscr{L}}
\newcommand{\mH}{\mathscr{H}}
\newcommand{\mS}{\mathscr{S}}

\newcommand{\Cbb}{\mathbb{C}}
\newcommand{\Nbb}{\mathbb{N}}
\newcommand{\Rbb}{\mathbb{R}}
\newcommand{\Sbb}{\mathbb{S}}
\newcommand{\tSbb}{\wtilde{\Sbb}}

\newcommand{\Zbb}{\mathbb{Z}}

\newcommand{\al}{\alpha}
\newcommand{\be}{\beta}
\newcommand{\la}{\lambda}
\newcommand{\Ups}{\Upsilon}
\newcommand{\bla}{\bar{\lambda}}
\newcommand{\de}{\delta}
\newcommand{\De}{\Delta}
\newcommand{\Ps}{\Psi}
\newcommand{\tPs}{\wtilde{\Ps}}
\newcommand{\hPs}{\what{\Ps}}
\newcommand{\Ph}{\Phi}
\newcommand{\vp}{\varphi}
\newcommand{\tvp}{\tilde\vp}
\newcommand{\hvp}{\hat\vp}
\newcommand{\Om}{\Omega}
\newcommand{\Th}{\Theta}

\newcommand{\ga}{\gamma}

\newcommand{\om}{\omega}

\newcommand{\stm}{\setminus}

\newcommand{\tZ}{\tilde{Z}}

\newcommand{\ql}{q(\la)}
\newcommand{\qp}{q_{\scriptscriptstyle+}}
\newcommand{\qm}{q_{\scriptscriptstyle-}}

\newcommand{\tf}{\tilde{f}}
\newcommand{\pzcf}{\mathpzc{f}}

\newcommand{\tg}{\tilde{g}}
\newcommand{\pzcg}{\mathpzc{g}}

\renewcommand{\th}{\tilde{h}}

\newcommand{\pzcm}{\mathpzc{m}}

\newcommand{\pzcT}{\mathpzc{T}}

\newcommand{\hu}{\hat{u}}
\newcommand{\pzcu}{\mathpzc{u}}

\newcommand{\tw}{\tilde{w}}
\newcommand{\hw}{\hat{w}}

\newcommand{\hx}{\hat{x}}

\newcommand{\ty}{\tilde{y}}
\newcommand{\hy}{\hat{y}}
\newcommand{\pzcy}{\mathpzc{y}}

\newcommand{\tz}{\tilde{z}}
\newcommand{\hz}{\hat{z}}

\newcommand{\pzcz}{\mathpzc{z}}

\newcommand{\ltp}{\ell^{\hspace{0.2mm}2}_{\Ps}}
\newcommand{\tltp}{\tilde{\ell}^{\hspace{0.3mm}2}_{\Ps}}
\newcommand{\tltpt}{\tilde{\ell}^{\hspace{0.3mm}2\times2}_{\Ps}}
\newcommand{\Tmax}{T_{\mathrm{max}}}
\newcommand{\Tmin}{T_{\mathrm{min}}}

\newcommand{\sZbb}{{\scriptscriptstyle{\Zbb}}}
\newcommand{\Iz}{\Ic_\sZbb}
\newcommand{\tIz}{\tIc_\sZbb}
\newcommand{\Izp}{\Ic_\sZbb^+}
\newcommand{\IzD}{\Ic_\sZbb^{\scriptscriptstyle{\rm D}}}

\newcommand{\mmatrix}[1]{\begin{pmatrix} #1
  \end{pmatrix}}
\newcommand{\msmatrix}[1]{\left(\begin{smallmatrix} #1
  \end{smallmatrix}\right)}  

\newcommand{\qtextq}[1]{\quad\text{#1}\quad}
\newcommand{\qtext}[1]{\quad\text{#1 }\ }

\DeclareMathOperator{\re}{Re}
\DeclareMathOperator{\im}{Im}
\newcommand{\rank}{\qopname\relax o{rank}}
\newcommand{\dom}{\qopname\relax o{dom}}
\newcommand{\cotan}{\qopname\relax o{cotan}}
\newcommand{\abs}[1]{|{#1}|}
\newcommand{\inner}[2]{\langle #1,#2 \rangle}
\newcommand{\innerP}[2]{\langle #1,#2 \rangle_{\Ps}}
\newcommand{\normP}[1]{\|{#1}\|_{\Ps}}
\newcommand{\normS}[1]{\|{#1}\|_2}
\newcommand{\normE}[1]{\|{#1}\|_2}

\newcommand{\Sla}[1]{\text{\rm(S$_{#1}$})}
\newcommand{\Slaf}[2]{\text{\rm(S$_{#1}^{#2}$)}}

\DeclareMathOperator{\adj}{adj}

\hypersetup{
pdfencoding=unicode,
pdftitle={\papertitle},
pdfauthor={\firstauthor{} \& \secondauthor},
pdfsubject={Department of Mathematics and Statistics, Faculty of Science, Masaryk University --- Department of Mathematics \& 
Statistics, Missouri University of Science and Technology},
pdfkeywords={\thekeywords},
pdfsubject=\paperabstract,
pdfencoding=auto
}

\begin{document}


\MakeFirstPageTwoAuthors

\section{Introduction}\label{S:intro}

This paper is devoted to the characterization of all self-adjoint extensions of the minimal linear relation associated with the 
discrete symplectic system
 \begin{equation*}\label{Sla}\tag{S$_\la$}
   z_{k}(\la)=\Sbb_{k}(\la)\,z_{k+1}(\la),\qquad \Sbb_{k}(\la)\coloneq\Sc_{k}+\la\Vc_{k},
 \end{equation*}
where $\la\in\Cbb$ is the spectral parameter, $\Sc_k$ and $\Vc_k$ are $2n\times 2n$ complex-valued matrices such that
 \begin{equation}\label{E:1.1}
  \Sc_k^*\Jc\Sc_k=\Jc,\quad \Vc_k^*\Jc\Sc_k\ \text{ is Hermitian},\quad \Vc_k^*\Jc\,\Vc_k=0
 \end{equation}
with the skew-symmetric $2n\times2n$ matrix $\Jc\!\coloneq\msmatrix{\phantom{-}0 & I\\ -I & 0}$ and the superscript $*$ 
denoting the conjugate transpose. Here $k$ belongs to a discrete interval $\Iz$, which is finite or unbounded from above. The 
conditions in \eqref{E:1.1} imply that $\Sbb_k(\la)$ satisfies the following symplectic-type identity
 \begin{equation}\label{E:1.2}
  \Sbb_k^*(\bla)\,\Jc\,\Sbb_k(\la)=\Jc,
 \end{equation}
which motivates the basic terminology for system \eqref{Sla}. Moreover, system \eqref{Sla} can be written as
 \begin{equation}\label{E:Sla.mL}
  \Jc\,[z_k(\la)-\Sc_k\,z_{k+1}(\la)]=\la\,\Ps_k\,z_k(\la),\quad \Ps_k\coloneq\Jc\Sc_k\,\Jc\,\Vc_k^*\,\Jc,
 \end{equation}
which gives to rise a linear map $\mL$ defined by the left-hand side of \eqref{E:Sla.mL}. Note that $\Ps$ also plays the 
role of the weight matrix in the associated semi-inner product (see Theorem~\ref{T:Lagrange}). Hence, we assume, in addition to 
\eqref{E:1.1}, that $\Ps_k$ is positive semidefinite on $\Iz$. 

System \eqref{Sla} is said to be in the ``time reversed'' form. Identity {\eqref{E:1.2}} implies that the matrix 
$\tSbb_k(\la)\coloneq\Sbb_k^{-1}(\la)=-\Jc\,\Sbb_k^*(\bla)\,\Jc$ exists for all $\la\in\Cbb$ and $k\in\Iz$, depends linearly on 
$\la$, and satisfies the same equality as in~\eqref{E:1.2}. Hence system \eqref{Sla} is equivalent with the (classical and more 
natural) ``forward'' discrete symplectic system
 \begin{equation*}\label{tSla}\tag{$\tilde{\rm{S}}_\la$}
  z_{k+1}(\la)=\tSbb_{k}(\la)\,z_{k}(\la),\qquad \tSbb_{k}(\la)\coloneq\tSc_{k}+\la\tVc_{k},
 \end{equation*}
i.e., $z(\la)$ solves \eqref{Sla} if and only if it solves system \eqref{tSla}. Moreover, the $\Ps$-norm of a~solution $z(\la)$ 
agrees with its $\tPs$-norm, where $\tPs_k$ denotes the weight matrix corresponding to system~\eqref{tSla} and $\tPs_k\geq0$ if 
and only if $\Ps_k\geq0$. This equivalence guarantees that the results of the Weyl--Titchmarsh theory 
for system~\eqref{tSla} established in \cite{slC.pZ10,rSH.pZ14:JDEA,mB.sS10} (see also \cite[Section~4]{rSH.pZ:MN}) are also true 
for system~\eqref{Sla}. 

Let us emphasize that analogously we could deal with system \eqref{tSla} instead of \eqref{Sla}. But the choice of system 
\eqref{Sla} is mainly motivated by the absence of the shift on the right-hand side of equality \eqref{E:Sla.mL} and in the 
associated semi-inner product, which produces more natural calculations, see \cite{slC.pZ15}. This is also the traditional 
approach in connection with the second order Sturm--Liouville difference equations (see, e.g. \cite{yS.hS11,dbH.rtL78}). On the 
other hand, systems~\eqref{Sla} and~\eqref{tSla} lead to different spaces of square summable sequences defined 
in~\eqref{E:ltp.def}.

Since the mapping associated with system \eqref{Sla} may be multivalued or non-densely defined, the approach dealing with 
linear relations instead of operators is utilized;  see \cite[Section~5]{slC.pZ15}. The study of linear relations associated with 
system~\eqref{Sla} begun in \cite{slC.pZ15} is continued in this paper with a characterization of self-adjoint extensions of 
symmetric linear relations. The description of self-adjoint extensions and their particular cases is a classic problem in the 
theory of differential and difference equations; see 
\cite{eaC.aD76,maE.dO.rpA11,xH.jS.aW.aZ12,xH.jS.aZ12,aW.jS.aZ09,yS.hS11,mbB.mB.hV14,oD.pH09,hdN.aZ92,hdN.aZ90,mM.aZ00,bmB.jsC05,
rsH.pZ10,maN68,jS86,hS.yS10,hS.yS11,dbH.amK.jkS87,dbH.anS97}. As in \cite{jS86,hS.yS11,gR.yS14,aW.jS.aZ09} our main result here 
is obtained by using square summable solutions of system~\eqref{Sla} and the Glazman--Krein--Naimark theory. 

In \cite{gR.yS14}, a characterization of self-adjoint extensions is given for linear Hamiltonian difference systems of the form
 \begin{equation}\label{E:1.2.hamilton}
  \De\mmatrix{x_k(\la)\\ u_k(\la)}=(H_k+\la W_k)\mmatrix{x_{k+1}(\la)\\ u_k(\la)},\quad
  H_k\coloneq\mmatrix{A_k & \phantom{-}B_k\\ C_k & -A^*_k},\quad
  W_k\coloneq\mmatrix{\phantom{-}0 & W_k^{[2]}\\ -W_k^{[1]} & 0},
 \end{equation}
where $B_k, C_k$, $W_k^{[1]}, W_k^{[2]}$ are $n\times n$ Hermitian matrices, $W_k^{[1]}\geq0$, $W_k^{[2]}\geq0$, and the matrix 
$I-A_k$ is invertible. We note that the underlying discrete interval considered in the latter reference can be also unbounded 
from 
below. An interesting overlap exists between the systems given in \eqref{Sla} and \eqref{E:1.2.hamilton}. System~\eqref{Sla} can 
be written as a linear Hamiltonian difference system only if the $n\times n$ matrix in the right-lower block of $\Sbb_k(\la)$ is 
invertible for all $\la\in\Cbb$ and $k\in\Iz$. However, in this instance the dependence on $\la$ may be nonlinear and the form of 
$W_k$ more general than in~\eqref{E:1.2.hamilton}. On the other hand, system \eqref{E:1.2.hamilton} can be written as \eqref{Sla} 
only if $W_k^{[2]}(I-A_k^*)^{-1}W_k^{[1]}\equiv0$. Without this additional assumption we obtain a discrete symplectic system 
with a special quadratic dependence on $\la$, see also \cite{rSH.pZ14:ICDEA,rSH.pZ15} for more details.

If we suppress the dependence on the spectral parameter, discrete symplectic systems, i.e., \eqref{Sla} or \eqref{tSla} with 
$\la=0$, represent the proper discrete counterpart of the linear Hamiltonian differential system (see, e.g. \cite{mB.oD97}). 
Hence system~\eqref{Sla} can be seen as a discrete analogue of the system 
 \begin{equation}\label{E:Hamiltion.R}
  z^\prime(t,\la)=\Jc\,[B(t)+\la A(t)]\,z(t,\la),
 \end{equation}
where $A(t)$, $B(t)$ are $2n\times 2n$  locally integrable, Hermitian matrix-valued functions (see 
Remark~\ref{R:canonical.transform}). But we point out the principal difference in the assumptions concerning the invertibility of 
the weight matrices $\Ps_k$ and $A(t)$. Hence we refer to \cite{rSH.pZ15}, where a connection between linear Hamiltonian 
differential and difference systems and discrete symplectic systems depending on the spectral parameter is discussed with using 
the time scale calculus, which provides suitable tools for this purpose.

The rest of the paper is organized as follows: In Section 2, we list notation used, introduce system~\eqref{Sla} precisely, and 
recall several results from the theory of linear relations. We also establish a~limit point criterion for system~\eqref{Sla} in 
Theorem~\ref{T:lpc}. In Section~\ref{S:main} we present the main result, Theorem~\ref{T:s-e.ext}, concerning the characterization 
of self-adjoint extensions of the minimal linear relation associated with system \eqref{Sla}. We apply this to a consideration of 
the $2\times 2$ (scalar) case for a finite discrete interval, and describe the Krein--von Neumann extension explicitly: see 
Theorems~\ref{T:5.2} and \ref{C:5.3}, and Example~\ref{Ex:krein.ext}. We note that there is no analogue of 
Theorems~\ref{T:lpc}, \ref{T:5.2}, \ref{C:5.3} and Example~\ref{Ex:krein.ext} in the setting of 
system~\eqref{E:1.2.hamilton}. Finally, Section~\ref{S:proof} is devoted to the proof of Theorem~\ref{T:s-e.ext}.

\section{Preliminaries}\label{S:prelim}

In the first part of this section we establish the basic notation. The real and imaginary parts of any $\la\in\Cbb$ are, 
respectively, denoted by $\re(\la)$ and $\im(\la)$, i.e., $\re(\la)\coloneq (\la+\bla)/2$ and 
$\im(\la)\coloneq (\la-\bla)/(2i)$. The symbols $\Cbb_+$ and $\Cbb_-$ mean, respectively, the upper and lower complex plane, 
i.e., $\Cbb_+\coloneq\{\la\in\Cbb\mid \im(\la)>0\}$ and $\Cbb_-\coloneq\{\la\in\Cbb\mid \im(\la)<0\}$. 

All matrices are considered over the field of complex numbers $\Cbb$. For $r,s\in\Nbb$ we denote by $\Cbb^{r\times s}$ the space 
of all complex-valued $r\times s$ matrices and $\Cbb^{r\times 1}$ will be abbreviated as $\Cbb^r$. For a given matrix 
$M\in\Cbb^{r\times s}$ we indicate by $M^\top$, $\overline{M}$, $M^*$, $\det M$, $\rank M$, $M\geq0$, $\adj(M)$, $\Rc(M)$, and 
$\dim\Rc(M)$, respectively, its transpose, conjugate, conjugate transpose, determinant, rank, positive definiteness, adjugate 
matrix, range (i.e., the space spanned by the columns of $M$) and the dimension of $\Rc(M)$. By $\normS{M}$, we denote the 
spectral norm for $M\in\Cbb^{n\times n}$,  i.e., $\normS{M}\coloneq\max\{\sqrt{\mu}\mid \mu\text{ is an eigenvalue of $M^*M$}\}$. 
This norm possesses the submultiplicative property, i.e., $\normS{MN}\leq\normS{M}\normS{N}$ for any $M,N\in\Cbb^{n\times n}$, 
and is the operator norm induced by the Euclidean norm on $\Cbb^n$, i.e., $\normE{v}=(v^*v)^{1/2}$ for any $v\in\Cbb^n$. 
Hence, we also have 
 \begin{equation}\label{E:norms}
  \normE{Mv}\leq\normS{M}\normE{v}
 \end{equation}
for any $M\in\Cbb^{n\times n}$ and $v\in\Cbb^{n}$. 

In addition, by $M_{u,v}$ we mean the submatrix of $M\in\Cbb^{r\times s}$ consisting of the first $u\leq r$ rows and of the first 
$v\leq s$ columns and we write only $M_u$ in the case $u=v$, i.e., for the $u$-th leading principal submatrix of $M$. The 
following relations are well known for any matrices $M\in\Cbb^{r\times s}$, $L\in\Cbb^{s\times p}$, and $Q\in\Cbb^{r\times q}$, 
 \begin{gather}
  \rank M+\rank L-s\leq \rank ML \leq \min\{\rank M, \rank L\}, \label{E:2.1}\\
  \rank M=\rank MM^*=\rank M^*M, \label{E:2.2}\\
  \rank (M,\ Q) +\dim[\Rc(M)\cap\Rc(Q)]=\rank M+\rank Q; \label{E:2.3}
 \end{gather}
e.g. \cite[Corollaries~2.5.1, 2.5.3, and 2.5.10 and Fact~2.11.9]{dsB09}.

Let $\mathcal{I}$ be an open or closed interval in $\Rbb$. Then, $\Iz\coloneq \mathcal{I}\cap\Zbb$ denotes the corresponding 
discrete interval. In particular, with $N\in\Nbb\cup\{0,\infty\}$, we shall be interested in discrete intervals of the form 
$\Iz\coloneq [0,N+1)_\Zbb$, in which case we define $\Izp\coloneq [0,N+1]_\sZbb$ with the understanding that $\Iz\equiv\Izp$ 
when $N=\infty$. Hence our system \eqref{Sla} will be considered on discrete intervals $\Iz$ which are finite or unbounded 
above. 

By $\Cbb(\Iz)^{r\times s}$ we denote the space of sequences defined on $\Iz$ of complex $r\times s$ matrices, where typically 
$r\in\{n,2n\}$ and $1\leq s\leq2n$. In particular, we write only $\Cbb(\Iz)^r$ in the case $s=1$. If $M\in\Cbb(\Iz)^{r\times s}$, 
then $M(k)\coloneq M_k$ for $k\in\Iz$ and if $M(\la)\in\Cbb(\Iz)^{r\times s}$, then $M(\la,k)\coloneq M_k(\la)$ for $k\in\Iz$ 
with 
$M_k^*(\la)\coloneq[M_k(\la)]^*$. If $M\in\Cbb(\Iz)^{r\times s}$ and $L\in\Cbb(\Iz)^{s\times p}$, then $MN\in\Cbb(\Iz)^{r\times 
p}$, where $(MN)_k\coloneq M_k N_k$ for $k\in\Iz$. The subspace of $\Cbb(\Iz)^{r\times s}$ consisting of all sequences compactly 
supported in $\Iz$ is denoted by $\Cbb_0(\Iz)^{r\times s}$. The forward difference operator acting on $\Cbb(\Iz)^{r\times s}$ is 
denoted by $\De$ where $(\De z)_k\coloneq\De z_k$. Finally, $z_k \big|_{m}^n\coloneq z_n - z_m$.

\subsection{Discrete symplectic systems}\label{Ss:systems}

In the previous section system~\eqref{Sla} was introduced through the matrices $\Sc,\Vc$ satisfying \eqref{E:1.1} and such that 
$\Ps$ given in \eqref{E:Sla.mL} is positive semidefinite. But according to \eqref{E:Sla.mL}, system~\eqref{Sla} can be determined 
also by $\Sc$ and a suitable matrix $\Ps$. This correspondence was shown in \cite[Subsection~2.1]{slC.pZ15} and it justifies 
the following hypothesis concerning the basic conditions for the coefficients of system~\eqref{Sla}. It guarantees that all the 
conditions in \eqref{E:1.1} are satisfied, which implies that any initial value problem associated with \eqref{Sla} is uniquely 
solvable on $\Iz$ for any initial value given at any $k_0\in\Izp$. This hypothesis is assumed throughout the paper.

\begin{hypothesis}\label{H:2.1}
 Let $n\in\Nbb$ and $\Iz$ be given. We have $\Sc,\Ps\in\Cbb(\Iz)^{2n\times2n}$ such that
  \begin{equation}\label{E:H.2.1}
   \Sc_k^*\Jc\Sc_k=\Jc,\quad \Ps_k^*=\Ps_k, \quad \Ps_k^*\,\Jc\,\Ps_k=0, \quad \Ps_k\geq0 \qtext{for all $k\in\Iz$.}
  \end{equation}
 Moreover, we define $\Sbb_{k}(\la)\coloneq\Sc_{k}+\la\Vc_{k}$ with $\Vc_k\coloneq -\Jc\,\Ps_k\,\Sc_k$ for all $k\in\Iz$.
\end{hypothesis}

Let us define the linear map 
 \begin{equation*}\label{E:mL}
  \mL:\Cbb(\Izp)^{2n}\to\Cbb(\Iz)^{2n},\qquad \mL(z)_k\coloneq\Jc(z_k-\Sc_k\,z_{k+1}).
 \end{equation*}
Then the nonhomogeneous problem
 \begin{equation}\label{Slaf}\tag{S$_\la^f$}
   z_k(\la)=\Sbb_k(\la)\,z_{k+1}(\la)-\Jc\,\Ps_k\,f_k, \quad k\in\Iz,
 \end{equation}
where $f\in\Cbb(\Iz)^{2n}$, can be written as 
 \begin{equation*}\label{E:2.5}
  \mL(z(\la))_k=\la\,\Ps_k\,z_k(\la)+\Ps_k\,f_k,\quad k\in\Iz;
 \end{equation*}
see \cite[Lemma~2.6]{slC.pZ15}. For convenience, we abbreviate $\mL^*(z)_k\coloneq [\mL(z)_k]^*$ and by \Slaf{\nu}{g} we will
refer to the nonhomogeneous system of the form \eqref{Slaf} with $\la$ replaced by $\nu$ and $f$ replaced by $g$. Analogous 
notation is employed also for system \eqref{Sla}, which corresponds to \Slaf{\la}{0}. We also suppress the dependence of $z(\la)$ 
on $\la$ when $\la=0$.

The following identity is crucial in the whole theory (see \cite[Theorem~2.5]{slC.pZ15} for its proof).

\begin{theorem}[Extended Lagrange identity]\label{T:Lagrange}
 Let $\la,\nu\in\Cbb$, $1\leq m\leq2n$, and $f,g\in\Cbb(\Iz)^{2n\times m}$. If $z(\la)\in\Cbb(\Izp)^{2n\times m}$ and 
 $u(\nu)\in\Cbb(\Izp)^{2n\times m}$ are solutions of systems~\eqref{Slaf} and \Slaf{\nu}{g}, respectively, then for any 
 $k,s,t\in\Iz$ such that $s\leq t$, we have
  \begin{align}
   \De[z_k^*(\la)\,\Jc u_k(\nu)]&=(\bla-\nu)\,z_k^*(\la)\,\Ps_k\,u_k(\nu)+f_k^*\,\Ps_k\,u_k(\nu)
                                                                                 -z_k^*(\la)\,\Ps_k\,g_k,\notag\\
   z_k^*(\la)\,\Jc u_k(\nu)\big|_{s}^{t+1}&=\sum_{k=s}^{t}\big\{(\bla-\nu)\,z_k^*(\la)\,\Ps_k\,u_k(\nu)+f_k^*\,\Ps_k\,u_k(\nu)
                                                                                 -z_k^*(\la)\,\Ps_k\,g_k\big\}.\label{E:2.7}
  \end{align}
 Especially, if $\nu=\bla$ and $f\equiv0\equiv g$, we get the Wronskian-type identity
  \begin{equation}\label{E:wronski.id}
   z_k^*(\la)\,\Jc\,u_k(\bla)=z_0^*(\la)\,\Jc\,u_0(\bla),\quad k\in\Izp.
  \end{equation}
\end{theorem}

Since we assume $\Ps_k\geq0$ on $\Iz$, Theorem~\ref{T:Lagrange} motivates the natural definition of the semi-inner product for 
$z,u\in\Cbb(\Izp)^{2n}$ as
 \begin{equation}\label{E:2.12}
  \innerP{z}{u}\coloneq\sum_{k\in\Iz} z_k^*\,\Ps_k\,u_k
 \end{equation}
and of the semi-norm $\normP{z}\coloneq\sqrt{\innerP{z}{z}}$. Then we denote by $\ltp$ the linear space of all square 
summable sequences defined on $\Izp$, i.e., 
 \begin{equation}\label{E:ltp.def}
  \ltp=\ltp(\Iz)\coloneq\{z\in\Cbb(\Izp)^{2n}\mid \normP{z}<\infty\}.
 \end{equation}

Identity \eqref{E:2.7} can be written as
 \begin{equation}\label{E:2.8}
  \big(z(\la),u(\nu)\big)_{\!k}\,\Big|_{s}^{t+1}=\sum_{k=s}^{t}\big\{\mL^*(z(\la))_k\,u_k(\nu)-z^*_k(\la)\,\mL(u(\nu))_k\big\},
 \end{equation}
where we use for any $z,u\in\Cbb(\Izp)^{2n}$ and $k\in\Izp$ the notation
 \begin{equation*}\label{E:2.10}
  (z,u)_k\coloneq z_k^*\,\Jc\,u_k.
 \end{equation*}
Moreover, under the assumptions of Theorem~\ref{T:Lagrange} with $\la=0=\nu$, $m=1$, $s=0$, and $t=N$ we get from \eqref{E:2.7} 
and \eqref{E:2.12} that
 \begin{equation}\label{E:2.13}
  (z,u)_{k}\,\Big|_{0}^{N+1}=\innerP{f}{u}-\innerP{z}{g},
 \end{equation}
where the left-hand side of \eqref{E:2.13} means $\lim_{k\to\infty} (z,u)_{k}-(z,u)_0$ if $\Iz=[0,\infty)_\sZbb$. Identity 
\eqref{E:2.13} shows that the latter limit exists finite whenever $z,u,f,g\in\ltp$.

\begin{remark}\label{R:canonical.transform}
 Similarly as in the continuous case, there exists a unitary map $Q:\Cbb(\Izp)^{2n}\to\Cbb(\Izp)^{2n}$ preserving the square 
 summability with respect to $\Ps$ and such that system~\Slaf{0}{f} can be written in the canonical form, i.e., with 
 $\Sc\equiv I$. Indeed, let $\Ph$ denote the fundamental matrix of system \Slaf{0}{0}  satisfying $\Ph_0=I$. Then, it is 
 invertible for all  $k\in\Izp$ with $\Ph_k^{-1}=-\Jc\,\Ph_k^*\,\Jc$ and this  inverse provides  the canonical transformation, 
 i.e., $Q=\Ph^{-1}$ with  $Q(z)_k\coloneq \Ph^{-1}_k z_k$. Hence system \Slaf{0}{f} is equivalent  with
  \begin{equation}\label{E:2.5A}
   -\Jc \De y_k=\hPs_k\,g_k,\quad k\in\Iz,
  \end{equation}
 where $y_k\coloneq Q(z)_k$, $g_k\coloneq Q(f)_k$, and $\hPs_k\coloneq \Ph_k^*\,\Ps_k\,\Ph_k$. One can easily verify that 
 $y\in\ell^2_{\hPs}$ if and only if $z\in\ltp$. System~\eqref{E:2.5A} can be seen as a discrete counterpart of the canonical 
 linear Hamiltonian differential system, i.e., nonhomogeneous system associated with \eqref{E:Hamiltion.R}, where 
 $B(t)\equiv0$; e.g. \cite[Subsection~2.2]{mL.mmM03} and the references therein.
\end{remark}

It is known that some Atkinson-type (or definiteness) condition is needed for the study of square summable solutions of discrete 
symplectic systems, see \cite{rSH.pZ14:JDEA,slC.pZ15}. These conditions guarantee that some (the ``weak'' condition) or all (the 
``strong'' condition) nontrivial solutions $z(\la)$ of \eqref{Sla} satisfy $\normP{z(\la)}\neq0$. The precise distinguishing 
between the weak and strong formulation of the Atkinson-type condition enables one to formulate some results of the 
Weyl--Titchmarsh theory for discrete symplectic systems with coupled (or jointly varying) endpoints, see \cite{rSH.pZ13}. On the 
other hand, the strong condition implies the equality between the number of linearly independent square summable solutions of 
system \eqref{Sla} and the deficiency index corresponding to the minimal linear relation associated with \eqref{Sla}, see 
\cite[Corollary~5.12]{slC.pZ15}. Since this relation shall be necessary for our treatment, we need the following hypothesis, see 
\cite[Section~3]{slC.pZ15}.

\begin{hypothesis}[Strong Atkinson condition]\label{H:SAC}
 There exists a finite interval $\IzD\coloneq[a,b]_\sZbb\subseteq\Iz$ such that for any $\la\in\Cbb$ every nontrivial solution 
 $z(\la)\in\Cbb(\Izp)^{2n}$ of system~\eqref{Sla} satisfies
  \begin{equation*}\label{E:2.14}
   \sum_{k=a}^b z^*_k(\la)\,\Ps_k\,z_k(\la)>0.
  \end{equation*}
\end{hypothesis}

The  positive semidefiniteness of $\Ps$ and Hypothesis~\ref{H:SAC} imply that
 \begin{equation}\label{E:2.15}
  0<\sum_{k=a}^b z^*_k(\la)\,\Ps_k\,z_k(\la)\leq \sum_{k\in \tIz} z^*_k(\la)\,\Ps_k\,z_k(\la)
 \end{equation}
for any discrete interval $\tIz$ such that $\IzD\subseteq\tIz\subseteq\Iz$. 

\begin{example}\label{Ex:2n.order}\ 
 \begin{enumerate}[leftmargin=10mm,topsep=1mm,label={{(\roman*)}}]
  \item As demonstrated in \cite[Example~3.4]{slC.pZ15}, the simplest example of system~\eqref{Sla} satisfying 
        Hypothesis~\ref{H:SAC} is represented by the scalar system 
         \begin{equation}\label{E:2nd.order=Sla}
          \mmatrix{x_k(\la)\\ u_k(\la)}=
          \mmatrix{1 & -1/p_{k+1}\\ -q_k+\la w_k & 1+(q_k-\la w_k)/p_{k+1}}\mmatrix{x_{k+1}(\la)\\ u_{k+1}(\la)}, \quad k\in\Iz,
         \end{equation}
        where $p_k,q_k,w_k$ are real-valued and such that $p_k\neq0$ on $\Izp$, $q_k$ is defined on $\Iz$, $w_k\geq0$ on $\Iz$, 
        and $w_k>0$ at least at two consecutive points of $\Iz$. In this case $\Ps_k=\msmatrix{w_k & 0\\ 0 & 0}$. System 
        \eqref{E:2nd.order=Sla} includes the second order Sturm--Liouville difference equation 
         \begin{equation}\label{E:2nd.S-L}
          -\De[p_k\,\De y_{k-1}(\la)]+q_k\,y_k(\la)=\la\,w_k\,y_k(\la),\quad k\in\Iz,
         \end{equation}
        (put $x_k=y_k$ and $u_k=p_k\,\De y_{k-1}$). Note that a~solution $y(\la)$ of the latter equation is defined on the 
        discrete interval $\{-1\}\cup\Izp$. 
  \item System~\eqref{E:2nd.order=Sla} is a particular case of system~\eqref{Sla} with the special linear dependence on $\la$, 
        i.e., 
         \begin{equation}\label{E:Sla.special}
          \text{\eqref{Sla}},\quad \Sc_k=\mmatrix{\Ac_k & \Bc_k\\ \Cc_k & \Dc_k},\quad 
          \Vc_k=\mmatrix{0 & 0\\ \Wc_k\,\Ac_k & \Wc_k\,\Bc_k},\quad k\in\Iz,
         \end{equation}
        where the $n\times n$ blocks are such that $\Sc_k$ satisfies the first equality in \eqref{E:H.2.1} 
        and $\Wc_k=\Wc^*_k\geq0$. Then Hypothesis~\ref{H:2.1} holds with $\Ps_k=\msmatrix{\Wc_k & 0\\ 0 & 0}$, because the first 
        equality in \eqref{E:H.2.1} equivalent with (suppressing the argument $k\in\Iz$)
         \begin{equation}\label{E:sympl.block.identity}
          \Ac^*\Dc-\Cc^*\Bc=I=\Ac\,\Dc^*-\Bc\,\Cc^*\qtextq{and}
          \Ac^*\Cc,\ \ \Bc^*\Dc,\ \ \Ac\,\Bc^*,\ \ \Cc\,\Dc^*\ \ \text{are Hermitian.}
         \end{equation}        
        In addition, if there exists an index $l\in\Iz\stm\{0\}$ such that the matrices
        $\Bc_{l-1},\Wc_{l-1},\Wc_{l}$ are invertible, then also Hypothesis~\ref{H:SAC} is satisfied, see 
        \cite[Theorem~3.11]{slC.pZ15}.
 \end{enumerate}
\end{example}

\begin{remark}\label{R:ql.def}
 If $\ql$ denotes the number of linearly independent square summable solution of system~\eqref{Sla} for $\la\in\Cbb$, i.e., 
  \begin{equation}\label{E:3.1}
   \ql\coloneq \dim Q(\la),\quad Q(\la)\coloneq\{z\in\ltp\mid\text{$z(\la)$ solves \eqref{Sla}}\},
  \end{equation}
 then under Hypothesis~\ref{H:SAC} (even its ``weak'' form) we have $n\leq\ql\leq 2n$ for all $\la\in\Cbb\stm\Rbb$, see 
 \cite[Section~4]{rSH.pZ14:JDEA} for more details. The geometrical background of this estimate leads to the classification of 
 system~\eqref{Sla} as being in the limit point case if $\ql=n$, and as being in the limit circle case if $\ql=2n$. Moreover, if 
 there exists $\la_0\in\Cbb$ such that $q(\la_0)=2n$, then $\ql\equiv2n$ on $\Cbb$, whether Hypothesis~\ref{H:SAC} is satisfied 
 or not, see \cite[Theorem~4.17]{rSH.pZ14:JDEA} and compare with the results in \cite{rSH.pZ:MN}. The latter statement is known 
 as the invariance of the limit circle case and a sufficient condition for this situation can be found in 
 \cite[Corollary~4.18]{rSH.pZ14:JDEA}. Consequently, under Hypothesis~\ref{H:SAC} (even its weak form) and with $n=1$, we obtain 
 the generalization of the well-known Weyl alternative: either all solutions of \eqref{Sla} belong to $\ltp$ for any 
 $\la\in\Cbb\stm\Rbb$, or there exists only one nontrivial solution in $\ltp$ for any $\la\in\Cbb\stm\Rbb$ (see 
 \cite[Corollary~4.19]{rSH.pZ14:JDEA}). Sufficient conditions for the invariance of $\ql$ in the case $q(\la_0)<2n$ remain open.
\end{remark}

The classical limit point criterion for linear Hamiltonian differential and difference systems \eqref{E:1.2.hamilton} and 
\eqref{E:Hamiltion.R} utilizes the minimal eigenvalue of the corresponding weight matrix. Unfortunately, similar criterion cannot 
be applied in the current setting, because the weight matrix $\Ps_k$ is always singular, see also 
\cite[Remark~4.16]{rSH.pZ14:JDEA}. In the following theorem we give conditions guaranteeing the invariance of the limit point 
case 
on $\Cbb\stm\Rbb$ for system~\eqref{Sla} with the special linear dependence on $\la$ as discussed in 
Example~\ref{Ex:2n.order}(ii). This statement is a discrete analogue of \cite[Theorem~5.6]{mL.mmM03}. 

\begin{theorem}\label{T:lpc}
 Let $\Iz=[0,\infty)_\sZbb$ and consider system~\eqref{E:Sla.special} such that $\Bc_k^*\,\Cc_k\equiv0$, $\Bc^*_k\,\Dc_k>0$, and 
 $\Wc_k>0$ for all $k\in\Iz$. If there exists $h\in\Cbb(\Iz)^1$ such that $h_k\geq h>0$ and 
  \begin{equation}\label{E:T.lpc.asmpt.1}
   \Ac_k^*\,\Cc_k\geq -h_k\,\Wc_{k+1},\quad
   \sum_{k=0}^{\infty} \frac{1}{g_k\sqrt{h_k}}=\infty,
  \end{equation}
 where $g_k\coloneq\max\big\{1,\big\|\Wc_{k+1}^{\,-1/2}(\Bc_k^*\,\Dc_k)^{-1/2}\big\|_2\big\}$, and a constant $T\geq0$ such 
 that 
  \begin{equation}\label{E:T.lpc.asmpt.2}
   \De\Big(\frac{1}{h_k}\Big)g_k\leq \frac{T}{\sqrt{h_k}}, \quad k\in\Iz,
  \end{equation}
 then system~\eqref{Sla} is in the limit point case for all $\la\in\Cbb\stm\Rbb$, i.e., $\ql=n$ for all $\la\in\Cbb\stm\Rbb$.
\end{theorem}
\begin{proof}
 The special structure of the coefficient matrices implies that system \eqref{Sla} can be written as
  \begin{align*}
   &x_{k}=\Ac_k\,x_{k+1}+\Bc_k\,u_{k+1},\\
   &u_{k}=(\Cc_k+\la\,\Wc_k\,\Ac_k)\,x_{k+1}+(\Bc_k+\la\,\Wc_k\,\Bc_k)u_{k+1}=\Cc_k\,x_{k+1}+\Bc_k\,u_{k+1}+\la\,\Wc_k\,x_{k},
  \end{align*}
 with $\Ps_k=\msmatrix{\Wc_k & 0\\ 0 & 0}$. The invertibility of $\Bc_k$ and $\Wc_k$ for all $k\in\Iz$ implies that 
 Hypothesis~\ref{H:SAC} holds, see \cite[Theorem~3.11]{slC.pZ15}. In accordance with \cite[Theorem~4.4]{rSH.pZ14:JDEA}, 
 and with $\ql$ defined in \eqref{E:3.1}, we have $\ql=n$ if and only if $\tZ_k(\la)\be\not\in\ltp$ for any 
 $\be\in\Cbb^n\stm\{0\}$, where $\tZ(\la)$ is the $2n\times n$ solution of system \eqref{Sla} determined by the initial condition 
 $\tZ_0(\la)=-\Jc\al^*$ with $\al\in\Cbb^{n\times2n}$ being such that $\al\,\al^*=I$ and  $\al\,\Jc\al^*=0$. Moreover, it is 
 sufficient to consider only $\la=\pm i$, because the number $\ql\geq n$ is constant in $\Cbb_+$ and $\Cbb_-$ by 
 \cite[Corollary~5.12]{slC.pZ15}. Hence, let $\be\in\Cbb^n\stm\{0\}$ and $\la\in\{\pm i\}$ be fixed. Let us denote
 $z_k\coloneq\msmatrix{x_k\\ u_k}=\tZ_k(\la)\,\be$ with the $n\times1$ components $x_k,u_k$ and $k\in\Iz$. Note that 
 $z_0^*\,\Jc z_0=0$. We show that under the current assumptions we have $z\not\in\ltp$. 
 
 Let us assume that $z\in\ltp$. By a direct calculation, we obtain from the block structure of the system and the identities in 
 \eqref{E:sympl.block.identity} that
  \begin{equation*}
   \De(x^*_k\,u_k)=-x^*_{k+1}\Ac_k^*\,\Cc_k\,x_{k+1}-x_{k+1}^*\Cc_k^*\,\Bc_k\,u_{k+1}-u_{k+1}^*\Bc_k^*\,\Cc_k\,x_{k+1}
                   -u_{k+1}^*\Bc_k^*\,\Dc_k\,u_{k+1}-\la\,x_k^*\,\Wc_k\,x_k.
  \end{equation*}
 Since $\Bc^*_k\,\Dc_k>0$ and $h_k>0$, the quantity 
 $\Fc_k(x,u)\coloneq \big(\sum_{j=0}^k \frac{1}{h_j}\,u^*_{j+1}\,\Bc_j^*\,\Dc_j\,u_{j+1}\big)^{\!1/2}\geq0$ is 
 well-defined. Then the latter equality and the assumption $\Bc_k^*\,\Cc_k\equiv0$ yield
  \begin{equation}\label{E:lpc.7}
   \Fc_k^2(x,u)
               =-\sum_{j=0}^k \frac{1}{h_j}\,x_{j+1}^*\,\Ac_j^*\,\Cc_j\,x_{j+1}
                -\la\sum_{j=0}^k \frac{1}{h_j}\,x_j^*\,\Wc_j\,x_j
                -\sum_{j=0}^k \frac{1}{h_j}\,\De(x^*_j\,u_j).
  \end{equation}
 From the Hermitian property and positive definiteness of $\Wc_k$ and $\Bc_k^*\,\Dc_k$, the Cauchy--Schwarz inequality, 
 inequality~\eqref{E:norms}, and the definition of $g_k$ we obtain 
  \begin{align}   
   \abs{x_{k+1}^*\,u_{k+1}}
    &=\abs{(\Wc_{k+1}^{1/2}\,x_{k+1})^*\,\Wc_{k+1}^{-1/2}\,(\Bc_k^*\,\Dc_k)^{-1/2}\,(\Bc_k^*\,\Dc_k)^{1/2}u_{k+1}}\notag\\
    &\leq\normE{\Wc^{1/2}_{k+1}\,x_{k+1}}\times
                                    \normE{\Wc^{-1/2}_{k+1}\,(\Bc_k^*\,\Dc_k)^{-1/2}\,(\Bc_k^*\,\Dc_k)^{1/2}u_{k+1}}\notag\\
    &\leq\normE{\Wc^{1/2}_{k+1}\,x_{k+1}}\times\normS{\Wc^{-1/2}_{k+1}\,(\Bc_k^*\,\Dc_k)^{-1/2}}
                                                                         \times\normE{(\Bc_k^*\,\Dc_k)^{1/2}u_{k+1}}\notag \\
    &\leq g_k\,\normE{\Wc^{1/2}_{k+1}\,x_{k+1}}\times\normE{(\Bc_k^*\,\Dc_k)^{1/2}u_{k+1}}.\label{E:lpc.8}
  \end{align}
 Hence the latter inequality, assumption~\eqref{E:T.lpc.asmpt.2}, the Cauchy--Schwarz inequality, and the inequality of 
 arithmetic  and geometric means $\sqrt{ab}\leq\frac{a+b}{2}$ yield
  \begin{align}
   \Bigg|\sum_{j=0}^k\! \De\Big(\frac{1}{h_j}\Big)\,x_{j+1}^*\,u_{j+1}\Bigg|
    &\leq \sum_{j=0}^k\! \De\Big(\frac{1}{h_j}\Big)\,\abs{x_{j+1}^*\,u_{j+1}}
    \leq \sum_{j=0}^k\! \De\Big(\frac{1}{h_j}\Big)\,g_j\,\normE{\Wc_{j+1}^{1/2}\,x_{j+1}}
                                                                        \times\normE{(\Bc_j^*\,\Dc_j)^{1/2}\,u_{j+1}}\notag\\
    &\leq \sum_{j=0}^k T\,\normE{\Wc_{j+1}^{1/2}\,x_{j+1}}\times h_j^{-1/2}\,\normE{(\Bc_j^*\,\Dc_j)^{1/2}\,u_{j+1}}\notag\\
    &\leq \bigg(T^2\sum_{j=0}^k \normE{\Wc_{j+1}^{1/2}\,x_{j+1}}^2\bigg)^{\!\!1/2}
                         \!\times\bigg(\sum_{j=0}^k h_j^{-1}\,\normE{(\Bc_j^*\,\Dc_j)^{1/2}\,u_{j+1}}^2\bigg)^{\!\!1/2}\notag\\
    &\leq \frac{1}{2}\bigg(T^2\sum_{j=0}^k \normE{\Wc_{j+1}^{1/2}\,x_{j+1}}^2+
                                   \sum_{j=0}^k h_j^{-1}\,\normE{(\Bc_j^*\,\Dc_j)^{1/2}\,u_{j+1}}^2\bigg)\notag\\  
    &\leq \frac{1}{2}\,\big(T^2\,\normP{z}^2+\Fc_k^2(x,u)\big).\label{E:lpc.9}
  \end{align}
 By using the summation by parts together with the inequalities $h_k\geq h$, \eqref{E:lpc.8}, and \eqref{E:lpc.9} we get
  \begin{align}
   \bigg|\re\sum_{j=0}^k \frac{1}{h_j}\,\De(x^*_j\,u_j)\bigg|
    &\leq\bigg|\sum_{j=0}^k \frac{1}{h_j}\,\De(x^*_j\,u_j)\bigg|
     \leq\bigg|\big[\,x^*_j\,u_j/h_j\big]_{0}^{k+1}-\sum_{j=0}^k \De\Big(\frac{1}{h_j}\Big)\,x^*_{j+1}\,u_{j+1} \bigg|\notag\\
    &\hspace*{-12mm}\leq\abs{x^*_0\,u_0/h_0}+\abs{x^*_{k+1}\,u_{k+1}/h_{k+1}}
                          +\bigg|\sum_{j=0}^k \De\Big(\frac{1}{h_j}\bigg)\,x^*_{j+1}\,u_{j+1}\Big|\notag\\
    &\hspace*{-12mm}\leq T_1+\frac{1}{h}\,g_k\,\normE{\Wc_{k+1}^{1/2}\,x_{k+1}}\times\normE{(\Bc_k^*\,\Dc_k)^{1/2}u_{k+1}}
                                                                 +\big(T^2\,\normP{z}^2+\Fc_k^2(x,u)\big)/2,\label{E:lpc.10}
  \end{align}
 where $T_1\coloneq |x^*_0\,u_0/h_0|$. Since 
 $\re\big(\Fc_k^2(x,u)\big)=-\sum_{j=0}^k \frac{1}{h_j}\,x_{j+1}^*\,\Ac_j^*\,\Cc_j\,x_{j+1}
 -\re\Big(\sum_{j=0}^k \frac{1}{h_j}\,\De(x^*_j\,u_j)\Big)$ and the inequality in \eqref{E:T.lpc.asmpt.1} implies 
 $-\sum_{j=0}^k \frac{1}{h_j}\,x_{j+1}^*\,\Ac_j^*\,\Cc_j\,x_{j+1}\leq \sum_{j=0}^k x_{j+1}^*\Wc_j\,x_{j+1}\leq \normP{z}^2$, it 
 follows from \eqref{E:lpc.7} and \eqref{E:lpc.10} that
  \begin{equation*}
   \frac{1}{2}\sum_{j=0}^k g_j^{-1}\,h_j^{-1/2}\,\Fc_j^2(x,u)\leq T_2\sum_{j=0}^k g_j^{-1}\,h_j^{-1/2}
     +\frac{1}{h}\sum_{j=0}^k h_j^{-1/2}\,\normE{\Wc_{j+1}^{1/2}\,x_{j+1}}\times\normE{(\Bc_j^*\,\Dc_j)^{1/2}u_{j+1}},
  \end{equation*}
 where $T_2\coloneq T_1+(1+T^2/2)\,\normP{z}^2$. Then with the aid of the Cauchy--Schwarz inequality we have
  \begin{align}
   G_k\coloneq \frac{1}{2}\sum_{j=0}^k g_j^{-1}\,h_j^{-1/2}\,[\Fc_j^2(x,u)-2\,T_2]
   &\leq\frac{1}{h}\,\bigg(\sum_{j=0}^k \normE{\Wc_{j+1}^{1/2}\,x_{j+1}}^2\bigg)^{\!\!1/2}\!\!
                \times\bigg(\sum_{j=0}^k \normE{(\Bc_j^*\,\Dc_j)^{1/2}u_{j+1}}^2\bigg)^{\!\!1/2}\notag\\
   &\leq\frac{1}{h}\,\normP{z}\,\Fc_k(x,u).\label{E:lpc.12}
  \end{align}
 In the next part we show that $\Fc_k^2(x,u)\leq2\,T_2$ for all $k\in\Iz$. Assume that there exists an index $m\in\Iz$ such 
 that $\Fc_m^2(x,u)>2\,T_2$. Since $\Fc_k^2(x,u)$ is nondecreasing, we have $\Fc_k^2(x,u)-2\,T_2>t$ for all 
 $k\in[m,\infty)_\sZbb$, where $t\coloneq \Fc_m^2(x,u)-2\,T_2$. Also $G_k$ is nondecreasing for all $k\in[m-1,\infty)_\sZbb$ 
 and for all $k\in[m,\infty)_\sZbb$ we obtain from \eqref{E:lpc.12} and the equality 
 $\Fc_k^2(x,u)=2\,g_k\,h_k^{1/2}\,\De G_{k-1}+2\,T_2$ that
  \begin{equation}\label{E:lpc.12A}
   h^2-2\,\normP{z}^2\,G_k^{-2}\,T_2 \leq 2\,G_k^{-2}\,\normP{z}^2\,g_k\,h_k^{1/2}\,\De G_{k-1}.
  \end{equation}
 In addition, $G_k\geq \frac{t}{2}\sum_{j=0}^k g_j^{-1}\,h_j^{-1/2}\to\infty$ for $k\to\infty$ by the second 
 part of \eqref{E:T.lpc.asmpt.1}. Now, let $0<a<2h^2$ be arbitrary and $l\in[m,\infty)_\sZbb$ be such that 
 $G_l\geq2\,\normP{z}\,T_2^{1/2}/\sqrt{2h^2-a}$. Then we have $a/2\leq h^2-2\,G_k^{-2}\,T_2\,\normP{z}^2$ for all 
 $k\in[l,\infty)_\sZbb$, which together with  \eqref{E:lpc.12A} yields for $k\in[l+1,\infty)_\sZbb$ that
  \begin{align*}
   \frac{a}{2}\,\sum_{j=l+1}^k\frac{1}{g_j\,h_j^{1/2}}
    &\leq \sum_{j=l+1}^k \frac{1}{g_j\,h_j^{1/2}} (h^2-2\,G_j^{-2}\,T_2\,\normP{z}^2)
     \leq \sum_{j=l+1}^k\! 2\,G_j^{-2}\,\normP{z}^2\,\De G_{j-1}\\
    &\leq 2\normP{z}^2\sum_{j=l+1}^k \frac{\De G_{j-1}}{G_j\,G_{j-1}}
     \leq -2\normP{z}^2\sum_{j=l+1}^k \De\bigg(\frac{1}{G_{j-1}}\bigg)
     \leq 2\,\normP{z}^2\,\frac{1}{G_l}<\infty.
  \end{align*}
 But it contradicts the second condition in \eqref{E:T.lpc.asmpt.1} for $k\to\infty$. Thus $\Fc_k^2(x,u)\leq2\,T_2$ for all 
 $k\in\Iz$, i.e.,
  \begin{equation}\label{E:lpc.13A}
   \sum_{j=0}^\infty h_j^{-1}\,u_{j+1}^*\,\Bc_j^*\,\Dc_j\,u_{j+1}\leq 2\,T_2<\infty.
  \end{equation} 
 Since system~\eqref{Sla} satisfies Hypothesis~\ref{H:SAC}, there exists $p\in\Iz$ such that 
 $\sum_{j=0}^{p} z_k^*\,\Ps_k\,z_k=T_3>0$. Hence the positive definiteness of $\Wc_k$ and the Lagrange identity in \eqref{E:2.7} 
 yield 
  \begin{equation}\label{E:lpc.14}
   \big|z_{k+1}^*\,\Jc z_{k+1}\big|=\Big|z_0^*\,\Jc z_0\pm 2i\sum_{j=0}^k z_j^*\,\Ps_j\,z_j\Big|
   =2\,\Big|\sum_{j=0}^k z_j^*\,\Ps_j\,z_j\Big|\geq2\,\Big|\sum_{j=0}^p z_j^*\,\Ps_j\,z_j\Big|=2\,T_3,
  \end{equation}
 for any $k\geq p$. Simultaneously, we get from \eqref{E:lpc.8} the estimate
  \begin{equation}\label{E:lpc.15}
   \big|z_{k+1}^*\,\Jc z_{k+1}\big|\leq 2\,\big|x^*_{k+1}\,u_{k+1}\big|
   \leq 2\,g_k\,h_k^{1/2}\,\normE{\Wc_{k+1}^{1/2}\,x_{k+1}}\times h_j^{-1/2}\,\normE{(\Bc_k^*\,\Dc_k)^{1/2}u_{k+1}}.
  \end{equation}
 The inequalities \eqref{E:lpc.13A}, \eqref{E:lpc.14}, \eqref{E:lpc.15}, and the Cauchy--Schwarz inequality imply for 
 $k\geq p$ that
  \begin{align*}
   \sum_{j=p}^k \frac{1}{g_j\,h_j^{1/2}}
   &\leq \sum_{j=p}^{k} \frac{2}{\big|z_{j+1}^*\,\Jc z_{j+1}\big|}\,\normE{\Wc_{k+1}^{1/2}\,x_{k+1}}
                                                                   \times h_j^{-1/2}\normE{(\Bc_k^*\,\Dc_k)^{1/2}u_{k+1}}\\
   &\leq \frac{1}{T_3} \sum_{j=p}^{k} \normE{\Wc_{k+1}^{1/2}\,x_{k+1}}
                                                                   \times h_j^{-1/2}\normE{(\Bc_k^*\,\Dc_k)^{1/2}u_{k+1}}\\
   &\leq \frac{1}{T_3} \bigg(\sum_{j=p}^{k} \normE{\Wc_{k+1}^{1/2}\,x_{k+1}}^2\bigg)^{\!\!1/2}
                             \!\!\times \bigg(\sum_{j=p}^{k} h_j^{-1}\normE{(\Bc_k^*\,\Dc_k)^{1/2}u_{k+1}}^2\bigg)^{\!\!1/2}
    \leq \frac{1}{T_3}\,\normP{z}\,\sqrt{2}\,T_2^{1/2}<\infty,   
  \end{align*}
 which (again) contradicts the second condition in~\eqref{E:T.lpc.asmpt.1} for $k\to\infty$. Hence $z\not\in\ltp$. 
 Since $\be$ and $\la$ were chosen arbitrarily, it follows that $\tZ(\la)\,\be\not\in\ltp$ for any $\be\in\Cbb^n\stm\{0\}$. 
 Therefore, system~\eqref{Sla} is in the limit point case for $\la\in\{\pm i\}$ and consequently for all $\la\in\Cbb\stm\Rbb$.
\end{proof}

Upon applying Theorem~\ref{T:lpc} to system~\eqref{E:2nd.order=Sla} with $q_k\equiv0$ we obtain the following corollary for a 
special case of the second order Sturm--Liouville difference equation~\eqref{E:2nd.S-L}, because one easily observes that  
$z(\la)\in\ltp$ if and only if $\sum_{k=0}^\infty \abs{y_k(\la)}^2\,w_k<\infty$, where 
$z_k(\la)=\big(y_k(\la), p_k\,\De y_{k-1}(\la)\big)^{\!\!\top}$ and $\Ps_k=\msmatrix{w_k & 0\\ 0 & 0}$.

\begin{corollary}\label{C:lpc}
 Let $\Iz=[0,\infty)_\sZbb$ and consider equation~\eqref{E:2nd.S-L} with $q_k\equiv0$, $p_k<0$ and $w_k>0$ for all $k\in\Iz$.
 If there exist $h_k\in\Cbb(\Iz)^1$ and a constant $T\geq0$ such that $h_k\geq h>0$ and 
  \begin{equation}\label{E:C.lpc.asmpt}
   \sum_{k=0}^{\infty} \frac{1}{g_k\,\sqrt{h_k}}=\infty,\quad
   \De\Big(\frac{1}{h_k}\Big)\,g_k\leq \frac{T}{\sqrt{h_k}}, \quad k\in\Iz,
  \end{equation}
 where $g_k\coloneq\max\big\{1,\big(-\frac{p_{k+1}}{w_{k+1}}\big)^{\!\scriptscriptstyle1/2}\big\}$, then 
 equation~\eqref{E:2nd.S-L} is in the limit point case for any $\la\in\Cbb\stm\Rbb$, i.e., there exists only one nontrivial 
 solution satisfying $\sum_{k=0}^\infty \abs{y_k(\la)}^2\,w_k<\infty$.
\end{corollary}

It was shown in \cite[Theorem~10]{dbH.rtL78}, see also \cite[Corollary~3.1]{hS.yS06:CMA}, that equation~\eqref{E:2nd.S-L} with 
$p_k\neq0$ and $w_k>0$ is in the limit point case for any $\la\in\Cbb\stm\Rbb$ if 
$\sum_{k=0}^\infty \frac{(w_k\,w_{k+1})^{1/2}}{\abs{p_{k+1}}}=\infty$. Corollary~\ref{C:lpc} partially generalizes this classical 
limit-point criterion as shown in the following example.

\begin{example}
 Let us consider the equation
  \begin{equation}\label{E:2nd.S-L.lpc}
   \text{\eqref{E:2nd.S-L}},\quad p_k\equiv-1,\quad q_k\equiv0,\quad w_k=1/(k+1)^2.
  \end{equation}
 Then the criterion from \cite[Theorem~10]{dbH.rtL78} cannot be applied, because
  \begin{equation*}
   \sum_{k=0}^\infty \frac{\big(w_k\,w_{k+1}\big)^{1/2}}{\abs{p_{k+1}}} 
   =\sum_{k=0}^\infty \sqrt{\frac{1}{(k+1)^2\,(k+2)^2}}=1<\infty.
  \end{equation*}
 On the other hand, the assumptions of Corollary~\ref{C:lpc} are satisfied with $h_k\equiv1$, $g_k=(k+2)$, and $T=0$, i.e., 
 equation~\eqref{E:2nd.S-L.lpc} is in the limit point case for all $\la\in\Cbb\stm\Rbb$. This fact can also be verified by 
 using the Weyl alternative; e.g. \cite[Theorem~5.6.1]{fvA64}. Indeed, equation~\eqref{E:2nd.S-L.lpc} with $\la=0$ has two 
 linearly independent solutions $y^{[1]}_k\equiv 1$ and $y^{[2]}_k=k$ for $k\in\{-1\}\cup\Iz$. Since  only $y^{[1]}$ is square 
 summable with respect to $w_k$, it follows from the Weyl alternative that equation~\eqref{E:2nd.S-L.lpc} has to be in the limit 
 point case for all $\la\in\Cbb\stm\Rbb$.
\end{example}

In the following lemma we establish a basic result concerning the solvability of a boundary value problem associated with 
\eqref{Sla}, which will be crucial in the proof of Lemma~\ref{L:2.6}. It provides the symplectic counterpart of the original 
Naimark's result known as the ``Patching lemma'', see \cite[Lemma~2 in Section~17.3]{maN68}. Analogous result for 
system~\eqref{E:1.2.hamilton} can be found in \cite[Lemma~3.3]{gR.yS14}.

\begin{lemma}\label{L:2.4}
 Let Hypothesis~\ref{H:SAC} be satisfied and a finite discrete interval $\tIz\coloneq[c,d]_\sZbb$ be given such that 
 $\IzD\subseteq\tIz\subseteq\Iz$. Then for any given $\al,\be\in\Cbb^{2n}$ there exists $f\in\Cbb(\tIz)^{2n}$ such that 
 the boundary value problem 
  \begin{equation}\label{E:2.16}
   \mL(z)_k=\Ps_k\,f_k,\quad z_c=\al,\quad z_{d+1}=\be,\quad k\in\tIz,
  \end{equation}
 has a solution $z\in\Cbb(\tIz^+)^{2n}$, where $\tIz^+\coloneq[c,d+1]_\sZbb$.
\end{lemma}
\begin{proof}
 Let $A$ be a $2n\times 2n$ matrix with the elements $a_{ij}\coloneq\sum_{k=c}^{d}\vp^{[i]*}_k\,\Ps_k\,\vp^{[j]}_k$ for 
 $i,j\in\{1,\dots,2n\}$, where $\vp^{[1]},\dots,\vp^{[2n]}\in\Cbb(\Izp)^{2n}$ are linearly independent solutions of 
 system \Sla{0}, i.e., $\mL(\vp^{[i]})_k=0$ for all $k\in\Iz$ and $i\in\{1,\dots,2n\}$. Then the homogeneous system of 
 algebraic equations $A\xi=0$, where $\xi=(\xi_1,\dots,\xi_{2n})^\top\in\Cbb^{2n}$, is equivalent with 
 $\sum_{k=c}^d \vp^*_k\,\Ps_k\,\vp_k=0$, where $\vp_k\coloneq\sum_{i=1}^{2n} \xi_i\,\vp_k^{[i]}$. Since $\vp$ also solves 
 system~\Sla{0}, it follows from Hypothesis~\ref{H:SAC} and inequality~\eqref{E:2.15} that $\vp$ is a trivial solution of 
 \Sla{0}, i.e., $\sum_{i=1}^{2n} \xi_i\,\vp_k^{[i]}\equiv0$, which implies that $\xi_i=0$ for all $i\in\{1,\dots,2n\}$. 
 It yields the invertibility of the matrix $A$.
 
 Hence there exists a unique solution $\eta=(\eta_1,\dots,\eta_{2n})^\top\in\Cbb^{2n}$ of the nonhomogeneous system
 of algebraic equations
  \begin{equation}\label{E:2.17}
   \eta^*A=\be^*\Jc\,\Phi_{d+1},
  \end{equation}
 where $\Ph\coloneq(\vp^{[1]*},\dots,\vp^{[2n]*})^*$ is a fundamental matrix of \Sla{0}. If we put $h_k^{[1]}\coloneq\Ph_k\,\eta$ 
 for $k\in\tIz$, we get from \eqref{E:2.17} for all $i\in\{1,\dots,2n\}$ that
  \begin{equation}\label{E:2.18}
   \sum_{k=c}^d h_k^{[1]*}\,\Ps_k\,\vp_k^{[i]}=\be^*\Jc\,\vp_{d+1}^{[i]}.
  \end{equation}
 Simultaneously Hypothesis~\ref{H:2.1} guarantees the existence of a unique solution $z^{[1]}\in\Cbb(\tIz^+)^{2n}$ of the 
 nonhomogeneous initial value problem
  \begin{equation*}\label{E:2.19}
   \mL(z^{[1]})_k=\Ps_k\,h_k^{[1]},\quad z^{[1]}_c=0,\quad k\in\tIz.
  \end{equation*}
 Then, for all $i\in\{1,\dots,2n\}$, the fact $\mL(\vp^{[i]})_k\equiv0$ and identity~\eqref{E:2.8} yield
  \begin{equation}\label{E:2.20}
   \sum_{k=c}^d h_k^{[1]*}\,\Ps_k\,\vp_k^{[i]}
    =\sum_{k=c}^d \big\{\mL^*(z^{[1]})_k\,\vp_k^{[i]}-z^{[1]*}_k\,\mL(\vp^{[i]})_k\big\}
    =(z^{[1]},\vp^{[i]})_k\big|_{c}^{d+1}
    =(z^{[1]},\vp^{[i]})_{d+1}.
  \end{equation}
 Upon combining \eqref{E:2.18} and \eqref{E:2.20} we obtain $z^{[1]}_{d+1}=\be$, which means that $z^{[1]}$ solves the boundary 
 value problem
  \begin{equation*}\label{E:2.21}
   \mL(z^{[1]})_k=\Ps_k\,h_k^{[1]},\quad z^{[1]}_c=0,\quad z^{[1]}_{d+1}=\be, \quad k\in\tIz.
  \end{equation*}
 Similarly, the nonhomogeneous system of algebraic equations $\om^*A=\al^*\Jc\,\Phi_c$ has a unique solution 
 $\om=(\om_1,\dots,\om_{2n})^\top\in\Cbb^{2n}$. Then with $h_k^{[2]}\coloneq\Ph_k\,\om$, $k\in\tIz$, we can calculate that 
 $z^{[2]}\in\Cbb(\tIz^+)^{2n}$, being the unique solution of
  \begin{equation*}\label{E:2.22}
   \mL(z^{[2]})_k=-\Ps_k\,h_k^{[2]},\quad z^{[2]}_{d+1}=0,\quad k\in\tIz,
  \end{equation*}
 also satisfies $z^{[2]}_c=\al$;  i.e., it solves the boundary value problem
  \begin{equation*}\label{E:2.23}
   \mL(z^{[2]})_k=-\Ps_k\,h_k^{[2]},\quad z^{[2]}_c=\al,\quad z^{[2]}_{d+1}=0, \quad k\in\tIz.
  \end{equation*}
 Thus, $z_k\coloneq z_k^{[1]}+z_k^{[2]}$, $k\in\tIz^+$, i.e., $z\in\Cbb(\tIz^+)^{2n}$, solves the boundary value problem 
 \eqref{E:2.16} with $f_k\coloneq h_k^{[1]}-h_k^{[2]}$ for $k\in\tIz$, i.e., $f\in\Cbb(\tIz)^{2n}$. 
\end{proof}

\subsection{Linear relations}\label{Ss:relations}

The theory of linear relations has been established as a suitable tool for the study of multi-valued or non-densely defined 
linear operators in a Hilbert space. Its history goes back to \cite{rA61} and the results were further developed e.g. in
\cite{sH.hsvdS.fhS09,eaC73:AMS,rC98,aD.hsvdS74}. In this subsection we recall the most relevant results from the theory of 
linear relations. A (closed) linear relation $\pzcT$ in a~Hilbert space $\mH$ over $\Cbb$ with the inner product 
$\inner{\cdot}{\cdot}$ is a (closed) linear subspace of the product space $\mH^2\coloneq \mH\times\mH$, i.e., the Hilbert 
space of all ordered pairs $\{\pzcz,\pzcf\}$ such that $\pzcz,\pzcf\in\mH$. By $\dom \pzcT$, $\ker \pzcT$, and $\overline{\pzcT}$ 
we mean, respectively, the domain of $\pzcT$, i.e., $\dom \pzcT\coloneq\{\pzcz\in\mH\mid \{\pzcz,\pzcf\}\in \pzcT\}$, the kernel 
of $\pzcT$, i.e., $\ker \pzcT\coloneq\{\pzcz\in\mH\mid \{\pzcz,0\}\in \pzcT\}$, and the closure of $\pzcT$. The sum 
$\pzcT+\mathpzc{U}$ and the algebraic sum $\pzcT\dotplus\mathpzc{U}$ are defined as
 \begin{gather*}
  \pzcT+\mathpzc{U}\coloneq\big\{\{\pzcz,\pzcf+\pzcg\}\mid\{\pzcz,\pzcf\}\in\pzcT,\ \{\pzcz,\pzcg\}\in\mathpzc{U} \big\},\\
  \pzcT\dotplus\mathpzc{U}\coloneq
     \big\{\{\pzcz+\pzcy,\pzcf+\pzcg\}\mid\{\pzcz,\pzcf\}\in\pzcT,\ \{\pzcy,\pzcg\}\in\mathpzc{U} \big\}.
 \end{gather*}

The adjoint $\pzcT^*$ of the linear relation $\pzcT$ is the closed linear relation defined by
 \begin{equation*}
  \pzcT^*\coloneq\big\{\{\pzcy,\pzcg\}\in\mH^2\mid 
                    \inner{\pzcz}{\pzcg}=\inner{\pzcf}{\pzcy}\ \text{for all }\ \{\pzcz,\pzcf\}\in \pzcT\big\}.
 \end{equation*}
A linear relation $\pzcT$ is said to be symmetric (or Hermitian) if $\pzcT\subseteq \pzcT^*$, and it is said to be self-adjoint 
if $\pzcT^*=\pzcT$. A symmetric linear relation $\pzcT_1$ is said to be a self-adjoint extension of $\pzcT$ if 
$\pzcT\subseteq \pzcT_1$ and $\pzcT_1^*=\pzcT_1$. 
For $\la\in\Cbb$ we define
 \begin{gather*}
  \pzcT-\la I\coloneq \big\{\{\pzcz,\pzcf-\la \pzcz\}\in\mH^2\mid \{\pzcz,\pzcf\}\in \pzcT\big\},\label{E:2.24a}\\
  M_{\la}(\pzcT)\coloneq \ker(\pzcT^*-\la I)=\{\pzcz\in\mH\mid \{\pzcz,\la \pzcz\}\in \pzcT^*\}.\label{E:2.24b}
 \end{gather*}
The number $d_\la(\pzcT)\coloneq \dim M_\la(\pzcT)$ is called the deficiency index of $\pzcT$ at $\la$ and the subspace
 \begin{equation*}\label{E:2.25}
  \Mc_\la(\pzcT)\coloneq \big\{\{\pzcz,\la \pzcz\}\in \pzcT^*\big\}
 \end{equation*}
denotes the defect space. It is known that the value of $d_\la(\pzcT)$ is constant in the upper and lower half plane of $\Cbb$, 
i.e., for $\la\in\Cbb_+$ and $\la\in\Cbb_-$. Hence we define the positive and negative deficiency indices as 
$d_\pm(\pzcT):=d_{\pm i}(\pzcT)$. If $\pzcT$ is a closed symmetric linear relation, then for every $\la\in\Cbb\stm\Rbb$ the 
following direct sum decomposition (a generalization of the von~Neumann formula)
 \begin{equation}\label{E:2.26}
  \pzcT^*=\pzcT\dotplus \Mc_\la(\pzcT) \dotplus \Mc_{\bla}(\pzcT)
 \end{equation}
holds, where the sum $\dotplus$ is orthogonal for $\la=\pm i$; e.g. \cite[Proposition~2.22]{mL.mmM03}. Moreover, for a~closed 
symmetric linear relation $\pzcT$ there is a self-adjoint extension if and only if $d_+(\pzcT)=d_-(\pzcT)$, see \cite[Corollary, 
pg.~34]{eaC73:AMS}.

The main results concerning the characterization of all self-adjoint extensions of the minimal linear relation associated with 
system \eqref{Sla} are obtained by applying the Glazman--Krein--Naimark theory for linear relations, which was established in 
\cite{yS12}. 

A complex linear space $\mS$ with a complex-valued function $[\,:\,]:\mS\times \mS\to\Cbb$ is called 
pre-symplectic if it possesses the conjugate bilinear and skew-Hermitian properties, i.e., for all $P,Q,R\in\mS$ and 
$\al\in\Cbb$ we have
 \begin{gather*}
  [P:Q+R]=[P:Q]+[P:R],\quad [P+Q:R]=[P:R]+[Q:R],\label{E:2.27}\\
  [\al P:Q]=\al\,[P:Q],\quad [P:\al Q]=\bar{\al}\,[P:Q],\label{E:2.28}\\
  [P:Q]=-\overline{[Q:P]};\label{E:2.29}
 \end{gather*}
see \cite{wnE.lM99:AMS} for more details. If we put $\mS=\mH^2$ and
 \begin{equation*}\label{E:2.30}
  [\{\pzcz,\pzcf\}:\{\pzcu,\pzcg\}]\coloneq \inner{\pzcf}{\pzcu}-\inner{\pzcz}{\pzcg}
 \end{equation*}
for $\{\pzcz,\pzcf\}$, $\{\pzcu,\pzcg\}\in\mH^2$, then $\mS$ and $[\,:\,]$ form the pre-symplectic space.

For a symmetric linear relation $\pzcT\subseteq\mH^2$ we have
 \begin{equation}\label{E:2.31}
  [\pzcT:\pzcT]=0=[\pzcT:\pzcT^*],\quad \overline{\pzcT}=\big\{\{z,f\}\in \pzcT^*\mid [\{z,f\}:\pzcT^*]=0\big\};
 \end{equation}
see \cite[Theorem~3.5]{yS12}. If, in addition, the linear relation $ \pzcT$ is closed and $d\coloneq d_+(\pzcT)=d_-(\pzcT)$, 
then the set $\{\be_j\}_{j=1}^d$ with $\be_j\in \pzcT^*$ for $j\in\{1,\dots,d\}$ such that
 \begin{enumerate}[leftmargin=10mm,topsep=1mm,parsep=1mm]
  \item $\be_1,\dots,\be_d$ are linearly independent in $\pzcT^*$ modulo $\pzcT$,
  \item $[\be_j:\be_i]=0$ for all $i,j\in\{1,\dots,d\}$,
 \end{enumerate}
is called GKN-set for the pair of linear relations $(\pzcT,\pzcT^*)$. The following theorem provides the necessary and sufficient 
conditions for a linear relation $\pzcT_1\subseteq\mH^2$ being a self-adjoint extension of $\pzcT$ (see 
\cite[Theorem~4.7]{yS12}).

\begin{theorem}\label{T:general.s-e.ext}
 Let $\pzcT\subseteq\mH^2$ be a closed symmetric linear relation such that $d_+(\pzcT)=d_-(\pzcT)=d$. A subspace 
 $\pzcT_1\subseteq\mH^2$ is a self-adjoint extension of $\pzcT$ if and only if there exists GKN-set $\{\be_j\}_{j=1}^d$ for 
 $(\pzcT,\pzcT^*)$ such that
  \begin{equation}\label{E:2.32}
   \pzcT_1=\{F\in \pzcT^*\mid [F:\be_j]=0\ \text{ for all } j=1,\dots,d\}.
  \end{equation}
\end{theorem}

A linear relation $\pzcT$ is called semibounded below, if there exists $a\in\Rbb$ such that
 \begin{equation}\label{E:semibound.def}
  \inner{\pzcz}{\pzcf}\geq a\,\inner{\pzcz}{\pzcz}\ \text{ for all }\{\pzcz,\pzcf\}\in\pzcT.
 \end{equation}
The number $\pzcm(\pzcT)\coloneq\sup\{a\in\Rbb\mid \text{\eqref{E:semibound.def} holds}\}$ is called the lower bound of 
$\pzcT$. If $\pzcm(\pzcT)>0$, the linear relation $\pzcT$ is said to be positive. Then, by analogy with the case of densely 
defined positive symmetric operators (see \cite[Theorem~5]{eaC.hsvdS78}),  the smallest and largest self-adjoint extensions of 
a positive symmetric linear relation are respectively known as the Krein--von~Neumann (or soft) extension $\pzcT_K$ and the 
Friedrichs (or hard) extension $\pzcT_F$. In particular, if $\pzcT$ is closed and $\pzcm(\pzcT)>0$, then the Krein--von~Neumann 
extension admits the representation
 \begin{equation}\label{E:krein.ext}
  \pzcT_K=\pzcT\dotplus(\ker \pzcT^*\times\{0\})
 \end{equation}
(see \cite[Corollary~1]{eaC.hsvdS78} and also \cite{sH.aS.hsvdS.hW07:JOT}).


\section{Main results}\label{S:main}

Since the weight matrix $\Ps$ is assumed to be only positive semidefinite in Hypothesis~\ref{H:2.1}, the space $\ltp$ is not a 
Hilbert space. Hence we need to consider the Hilbert space of equivalence classes. It is the quotient space obtained by factoring 
out the kernel of the semi-norm $\normP{\cdot}$, i.e., the space
 \begin{equation*}\label{E:2.33}
  \tltp=\tltp(\Iz)\coloneq \ltp\big/\big\{z\in\Cbb(\Izp)^{2n}\mid \ \normP{z}=0\big\}
 \end{equation*}
with the inner product $\innerP{\tz}{\tf}\coloneq\innerP{z}{f}$, where $z$ and $f$ are elements of the equivalence classes $\tz$, 
$\tf\in\tltp$. Note that the value $\normP{\tz}$ for $z\in\Cbb(\Izp)^{2n}$ does not depend on $z_{N+1}$ in the case of $\Iz$ 
being 
a~finite discrete interval, which implies that the sequences $z,y\in\Cbb(\Izp)^{2n}$ such that $z_k\neq y_k$ only for $k=N+1$, 
belong to the same equivalence class. We also introduce the space $\ell^2_{\Ps,0}$ as
 \begin{equation*}
  \ell^2_{\Ps,0}\coloneq
   \begin{cases}
    \big\{z\in \Cbb_0(\Izp)^{2n}\mid  z_0=0\big\} & \text{if } N=\infty,\\[1mm]
    \big\{z\in \Cbb_0(\Izp)^{2n}\mid  z_0=0,\ z_{N+1}=0 \big\} & \text{if } N\in\Nbb\cup\{0\}.
   \end{cases}
 \end{equation*}
Moreover, the corresponding function $[\,:\,]:\tltpt\times\tltpt\to\Cbb$ for the pre-symplectic space associated with  
$\tltpt\coloneq\tltp\times \tltp$ is given by
 \begin{equation*}\label{E:2.37}
  [\{\tz,\tf\}:\{\tw,\tg\}]\coloneq\innerP{\tf}{\tw}-\innerP{\tz}{\tg}.
 \end{equation*}

Linear relations associated with system \eqref{Sla} were introduced and studied in \cite[Section~5]{slC.pZ15}. The maximal 
linear relation in $\tltpt$ is defined as
 \begin{equation*}\label{E:2.34}
  \Tmax\coloneq \big\{\{\tz,\tf\}\in\tltpt\mid \text{there exists $u\in\tz$ such that } 
                                                     \mL(u)_k=\Ps_k\,f_k\ \text{ for all } k\in\Iz\big\}.
 \end{equation*}
Observe that the above definition does not depend on the particular choice of $f\in\tf$. Moreover, Hypothesis~\ref{H:SAC} is 
satisfied if and only if for any $\{\tz,\tf\}\in\Tmax$ there exists unique $u\in\tz$ such that $\mL(u)_k=\Ps_k\,f_k$ for all 
$k\in\Iz$, see \cite[Theorem~5.2]{slC.pZ15}. Henceforth, this unique element shall be denoted as $\hz$. Then 
identity~\eqref{E:2.13} yields for any $\{\tz,\tf\},\{\tw,\tg\}\in\Tmax$ that
 \begin{equation}\label{E:2.38}
  [\{\tz,\tf\}:\{\tw,\tg\}]=\innerP{f}{\hw}-\innerP{\hz}{g}=(\hz,\hw)_k\big|_{0}^{N+1},
 \end{equation}
where $f\in\tf$ and $g\in\tg$ are arbitrary representatives. Thus, under Hypothesis~\ref{H:SAC}, we obtain from Lemma~\ref{L:2.4} 
the following statement, compare with \cite[Remark~3.2]{gR.yS14} and 
\cite[Lemma~3.3]{yS.hS11}.

\begin{lemma}\label{L:2.6}
 Let Hypothesis~\ref{H:SAC} be satisfied. Then for any pairs $\{\tz,\tf\},\{\tw,\tg\}\in\Tmax$ there exists $\{\ty,\th\}\in\Tmax$ 
 such that
  \begin{equation*}
   \hy_k=\begin{cases}
        \hz_k, & k\in[0,c]_\sZbb\cap\Iz,\\
        \hw_k, & k\in[d+1,\infty)_\sZbb\cap\Izp,
       \end{cases}
  \end{equation*}
 where $c\in[0,a]_\sZbb\cap\Iz$, $d\in[b,\infty)_\sZbb\cap\Izp$ with $a,b$ determining the interval $\IzD$ in 
 Hypothesis~\ref{H:SAC}.
 
 In particular, for $i\in\{1,\dots,2n\}$ there exists $\{\tz^{[i]},\tf^{[i]}\}\in\Tmax$ such that $\hz^{[i]}_0=e_i$ and 
 $\hz_k^{[i]}=0$ for $k\in[d+1,\infty)_\sZbb\cap\Izp$, where $e_i=(0,\dots,1,\dots,0)^\top\in\Cbb^{2n}$ is the $i$-th canonical 
 unit vector. If, in addition, $N\in\Nbb\cup\{0\}$, i.e., $\Iz$ is a finite discrete interval, then there exists 
 $\{\ty^{[i]},\th\}\in\Tmax$ such that $\hy^{[i]}_{N+1}=e_i$ and $\hy_k^{[i]}=0$ for $k\in[0,c]_\sZbb\cap\Iz$.
\end{lemma}
\begin{proof}
 Let $\tIz$ be a finite discrete interval as in Lemma~\ref{L:2.4}, the pairs $\{\tz,\tf\},\{\tw,\tg\}\in\Tmax$ be arbitrary, and 
 define $\al\coloneq\hz_c$, $\be\coloneq \hw_{d+1}$. Then, by the latter lemma there exist sequences $l\in\Cbb(\tIz)^{2n}$ and  
 $v\in\Cbb(\tIz^+)^{2n}$ such that
  \begin{equation*}
   \mL(v)_k=\Ps_k\,l_k,\quad v_c=\al,\quad v_{d+1}=\be,\quad k\in\tIz.
  \end{equation*}
 Putting
  \begin{equation*}
   y_k\coloneq
    \begin{cases}
     \hz_k, & k\in[0,c]_\sZbb\cap\Iz,\\
     v_k, & k\in[c+1,d]_\sZbb\cap\Iz,\\
     \hw_k, & k\in[d+1,\infty)_\sZbb\cap\Izp
    \end{cases}\qquad
   h_k\coloneq
    \begin{cases}
     f_k, & k\in[0,c-1]_\sZbb\cap\Iz,\\
     l_k, & k\in[c,d]_\sZbb\cap\Iz,\\
     g_k, & k\in[d+1,\infty)_\sZbb\cap\Iz,
    \end{cases}
  \end{equation*}
 it can be verified by a direct calculation that $y,g\in\ltp$ and that they satisfy $\mL(y)_k=\Ps_k\,h_k$ for $k\in\Iz$, i.e., 
 $\{\ty,\th\}\in\Tmax$ with $\hy_k\equiv y_k$. The second part of the statement follows directly from Lemma~\ref{L:2.4}.
\end{proof}

The minimal linear relation is defined as $\Tmin\coloneq \overline{T_0}$, where $T_0$ is the pre-minimal linear relation
 \begin{equation*}\label{E:2.35}
  T_0\coloneq \big\{\{\tz,\tf\}\in\tltpt\mid \text{there exists $u\in\tz\cap\ell^2_{\Ps,0}$ such that } 
                                                             \mL(u)_k=\Ps_k\,f_k\ \text{ for all } k\in\Iz\big\}.
 \end{equation*}
It was shown in \cite[Theorem~5.10]{slC.pZ15} that
 \begin{equation}\label{E:2.36}
  T_0^*=\Tmin^*=\Tmax,
 \end{equation}
which implies that $\Tmin$ is a closed and symmetric linear relation. Moreover, the following theorem provides a more explicit
characterization of $\Tmin$; cf. \cite[Theorem~3.2]{gR.yS14}.

\begin{theorem}\label{T:2.7}
 Let Hypothesis~\ref{H:SAC} be satisfied. Then,
  \begin{equation}\label{E:2.39}
   \Tmin=\big\{\{\tz,\tf\}\in\Tmax\mid \hz_0=0=(\hz,\hw)_{N+1}\ \text{ for all } \tw\in\dom\Tmax\big\},
  \end{equation}
 which in the case of $\,\Iz$ being a finite discrete interval reduces to
  \begin{equation}\label{E:2.40}
   \Tmin=\big\{\{\tz,\tf\}\in\Tmax\mid \hz_0=0=\hz_{N+1}\big\}.
  \end{equation}
\end{theorem}
\begin{proof}
 Since $\Tmin=\overline{(\Tmin)}$ by the definition, identities~\eqref{E:2.31}, \eqref{E:2.38}, and \eqref{E:2.36} yield
  \begin{equation}\label{E:Tmin.eq.pf}
   \Tmin=\big\{\{\tz,\tf\}\in\Tmax\mid (\hz,\hw)_k\big|_{0}^{N+1}=0\ \text{ for all } \hw\in\dom\Tmax\big\}.
  \end{equation}
 Let $T$ be the linear relation on the right-hand side of \eqref{E:2.39}. Then, it is obvious that $T\subseteq\Tmin$. On the 
 other hand, let $\{\tz,\tf\}\in\Tmin$ be fixed. Then, $(\hz,\hw)_k\big|_{0}^{N+1}=0$ for all $\hw\in\dom\Tmax$ by 
 \eqref{E:Tmin.eq.pf}. By Lemma~\ref{L:2.6}, for any $\{\tw,\tg\}\in\Tmax$ there exists $\{\ty,\th\}\in\Tmax$ such that $\hy_k=0$ 
 for $k\in[0,c]_\sZbb\cap\Iz$ and $\hy_k=\hw_k$  for $k\in[d+1,\infty)_\sZbb\cap\Izp$. Hence $(\hz,\hw)_0=(\hz,\hw)_{N+1}=0$ for 
 all $\hw\in\dom\Tmax$. From the second part of Lemma~\ref{L:2.6} we get $\hz_0=0$, because there exists 
 $\{\tz^{[i]},\tf^{[i]}\}\in\Tmax$ such that $\hz^{[i]}_0=e_i$. Therefore, $T=\Tmin$. If, in addition, $\Iz$ is a finite 
 discrete interval, i.e., $N\in\Nbb$, then $\dom\Tmax$ contains also $\ty$ such that $\hy_{N+1}=e^{[i]}$, $i\in\{1,\dots,2n\}$,  
 by the last part of Lemma~\ref{L:2.6}. Hence equality \eqref{E:2.40} holds.
\end{proof}

By \cite[Corollary~5.12]{slC.pZ15}, Hypothesis~\ref{H:SAC} is equivalent with the equality $\ql=d_\la(\Tmin)$, which means 
that the number of the linearly independent square summable solutions of \eqref{Sla} is constant in $\Cbb_+$ and $\Cbb_-$. 
Therefore the numbers $\qp\coloneq \ql$ for $\la\in\Cbb_+$ and $\qm\coloneq \ql$ for $\la\in\Cbb_-$ are well-defined for $q(\la)$ 
given in \eqref{E:3.1}. Let $\la_0\in\Cbb_+$ be fixed. Then system \Sla{\la_0} has $\qp$ linearly independent square summable 
solutions, which we denote as $v^{[1]}(\la_0),\dots,v^{[\qp]}(\la_0)$, and similarly system~\Sla{\bla_0} has $\qm$ linearly 
independent square summable solutions, which we denote as $w^{[1]}(\bla_0),\dots,w^{[\qm]}(\bla_0)$. Let
 \begin{gather}\label{E:3.3+3.4}
  \vp^{[i]}_k\coloneq v_k^{[i]}(\la_0),\quad \vp^{[j+\qp]}_k\coloneq w_k^{[j]}(\bla_0),\quad i=1,\dots,\qp,\quad 
  j=1,\dots,\qm,\quad k\in\Izp,
 \end{gather}
and $\Ph_k\coloneq(\Ph_k^+,\Ph_k^-)\in\Cbb(\Izp)^{2n\times p}$, where $\Ph_k^+\coloneq(\vp^{[1]}_k,\dots,\vp^{[\qp]}_k)$ and 
$\Ph_k^-\coloneq(\vp^{[1+\qp]}_k,\dots,\vp^{[p]}_k)$ with $2n\leq p\coloneq \qp+\qm\leq 4n$. Then for $i\in\{1,\dots,\qp\}$ and 
$j\in\{\qp+1,\dots,p\}$ we have $\{\tvp^{[i]},\la_0\,\tvp^{[i]}\}\in\Tmax$ and $\{\tvp^{[j]},\bla_0\,\tvp^{[j]}\}\in\Tmax$ with 
$\hvp^{[l]}\equiv\vp^{[l]}$ for $l=1,\dots,p$. We also define the matrix
 \begin{equation}\label{E:Om.def}
  \Om=\mmatrix{\Om^{[1,1]} & \Om^{[1,2]}\\ \Om^{[2,1]} & \Om^{[2,2]}}
      \coloneq\mmatrix{(\vp^{[1]},\vp^{[1]})_{N+1} & \hdots & (\vp^{[1]},\vp^{[p]})_{N+1}\\
               \vdots & \ddots & \vdots\\
              (\vp^{[p]},\vp^{[1]})_{N+1} & \hdots & (\vp^{[p]},\vp^{[p]})_{N+1}}\in\Cbb^{p\times p},
 \end{equation}
where $\Om^{[1,2]}\in\Cbb^{\qp\times\qm}$. Note that the elements $\om_{ij}\coloneq (\vp^{[i]},\vp^{[j]})_{N+1}$ exist finite 
for all $i,j=1,\dots,p$ by identity~\eqref{E:2.13}. Moreover, from \eqref{E:wronski.id} one easily concludes that the matrix 
$\Om^{[1,2]}$ consists of the elements $(\vp^{[i]},\vp^{[j]})_{N+1}=(\vp^{[i]},\vp^{[j]})_0$ for $i\in\{1,\dots,\qp\}$ and 
$j\in\{\qp+1,\dots,p\}$.

Upon combining \eqref{E:2.36} and \eqref{E:2.26} we get that any $\{\tz,\tf\}\in\Tmax$ can be written as
 \begin{equation}\label{E:3.5}
  \hz_k=\hy_k+\sum_{j=1}^{p} \xi_j \vp^{[j]}_k,\quad k\in\Izp,
 \end{equation}
where $\hy\in\dom\Tmin$ and $\xi_1,\dots,\xi_p\in\Cbb$ are determined uniquely. Especially, for 
$\{\tz^{[i]},\tf^{[i]}\}\in\Tmax$ (see Lemma~\ref{L:2.6}), we get the unique expression
 \begin{equation}\label{E:3.6}
  \hz^{[i]}_k=\hy^{[i]}_k+\sum_{j=1}^{p} \xi_{i,j} \vp^{[j]}_k,\quad k\in\Izp,\quad i=1,\dots,2n.
 \end{equation}
If we put $Z_k\coloneq(\hz^{[1]}_k,\dots,\hz^{[2n]}_k)$ for $k\in\Izp$, then identity~\eqref{E:3.6} implies
 \begin{equation}\label{E:3.6A}
  Z_k=Y_k+\Ph_k\,\Xi^\top,
 \end{equation}
where $Y_k\coloneq(\hy^{[1]}_k,\dots,\hy^{[2n]}_k)\in\Cbb(\Izp)^{2n\times2n}$ and the matrix $\Xi\in\Cbb^{2n\times p}$ consists 
of the elements $\xi_{i,j}$. In particular, for $k=0$ we obtain $I=Y_0+\Ph_0\, \Xi^\top$, which together with \eqref{E:2.39} 
yields $I=\Ph_0\,\Xi^\top$, i.e., $\rank\Xi=2n$ by the second inequality in \eqref{E:2.1}. From the definition of $\hz^{[i]}$, 
its expression in \eqref{E:3.6}, and identity \eqref{E:2.39} we have
 \begin{equation*}\label{E:3.9+10}
  0=(\hz^{[i]},\vp^{[l]})_{N+1}=(\hy^{[i]},\vp^{[l]})_{N+1}+\sum_{j=1}^p \overline{\xi_{i,j}}\,(\vp^{[j]},\vp^{[l]})_{N+1}
   =\sum_{j=1}^p \overline{\xi_{i,j}}\,(\vp^{[j]},\vp^{[l]})_{N+1}
 \end{equation*}
for all $i\in\{1,\dots,2n\}$ and any $l\in\{1,\dots,p\}$, i.e., $\overline{\Xi}\,\Om=0$. Since $\rank\Xi=2n$, the first 
inequality in 
\eqref{E:2.1} implies
 \begin{equation*}\label{E:3.13}
  \rank\Om\leq p-2n.
 \end{equation*}
On the other hand, the equality $\Om^{[1,2]}=\Ph_0^{+*}\Jc\,\Ph_0^-$ and the first inequality in \eqref{E:2.1} yield
 \begin{equation*}\label{E:3.15}
  \rank\Om^{[1,2]}\geq p-2n.
 \end{equation*}
Therefore, $\rank\Om=p-2n=\rank\Om^{[1,2]}$. Since $p-2n\leq\qp$ and $p-2n\leq\qm$, we may assume, without loss of generality, 
that $\vp^{[1]},\dots,\vp^{[\qp]}$ are arranged such that
 \begin{equation}\label{E:rank.Om.1.2}
  \rank\Om^{[1,2]}_{p-2n,\,\qm}=p-2n.
 \end{equation}

The main result concerning the characterization of all self-adjoint extension of $\Tmin$ is stated in the following theorem and 
its proof is given in Section~\ref{S:proof}; cf. \cite[Theorem~5.7]{gR.yS14}. Recall that for the existence of a 
self-adjoint extension it is essential to assume $\qp=\qm$.

\begin{theorem}\label{T:s-e.ext}
 Let Hypothesis~\ref{H:SAC} be satisfied, equality $\qp=\qm\eqcolon q$ hold and assume that the solutions 
 $\vp^{[1]},\dots,\vp^{[q]}$ are arranged such that \eqref{E:rank.Om.1.2} holds. Then a linear relation $T\subseteq \tltpt$ is 
 a~self-adjoint extension of $\Tmin$ if and only if there exist matrices $M\in\Cbb^{q\times 2n}$ and $L\in\Cbb^{q\times (2q-2n)}$ 
 such that
  \begin{equation}\label{E:4.1}
   \rank(M,L)=q,\quad M\Jc M^*-L\,\Om_{2q-2n}\,L^*=0,
  \end{equation}
 and 
  \begin{equation}\label{E:4.2}
   T=\Bigg\{\{\tz,\tf\}\in\Tmax\mid M \hz_0-L\msmatrix{(\vp^{[1]},\hz)_{N+1}\\ \vdots\\ (\vp^{[2q-2n]},\hz)_{N+1}}=0\Bigg\}.
  \end{equation}
\end{theorem}

\begin{remark}\label{R:Ups.def}
 If, in addition to the assumptions of Theorem~\ref{T:s-e.ext}, there exists $\nu\in\Rbb$ such that \Sla{\nu} has $q$ linearly 
 independent square summable solutions (suppressing the argument $\nu$) $\Th^{[1]},\dots,\Th^{[q]}$, then the statement of 
 Theorem~\ref{T:s-e.ext} can be formulated by using these solutions, which are (without loss of generality) arranged such that 
 the submatrix $\Ups_{2q-2n}$ has the full rank, where 
  \begin{equation*}
   \Ups\coloneq\mmatrix{(\Th^{[1]},\Th^{[1]})_{N+1} & \hdots & (\Th^{[1]},\Th^{[q]})_{N+1}\\
               \vdots & \ddots & \vdots\\
               (\Th^{[q]},\Th^{[1]})_{N+1} & \hdots & (\Th^{[q]},\Th^{[q]})_{N+1}},
  \end{equation*}
 see Lemma~\ref{L:3.4}. Moreover, the Wronskian-type identity \eqref{E:wronski.id} yields that $\Ups=\Th^*_0\,\Jc\,\Th_0$, where 
 $\Th_k\coloneq(\Th_k^{[1]},\dots,\Th_k^{[q]})$ for $k\in\Izp$.
\end{remark}
 
In the next part we discuss several special cases of Theorem~\ref{T:s-e.ext}. If system \eqref{Sla} is in the limit point case 
for all $\la\in\Cbb\stm\Rbb$, i.e., $\qp=\qm=n$, then the boundary conditions at $N+1$ (which is necessary equal to $\infty$) are 
superfluous as stated in the following corollary; cf. \cite[Theorem~5.9]{gR.yS14}. \emph{This situation occurs, e.g., when 
the assumptions of Theorem~\ref{T:lpc} are satisfied}. The proof follows directly from Theorem~\ref{T:s-e.ext}.

\begin{corollary}\label{C:s-e.LPC}
 Let Hypothesis~\ref{H:SAC} be satisfied and $\qp=\qm=n$ hold. Then a linear relation $T\subseteq \tltpt$ 
 is a self-adjoint extension of $\Tmin$ if and only if there exists a matrix $M\in\Cbb^{n\times 2n}$ such that
  \begin{equation*}\label{E:4.3}
   \rank M=n,\quad M \Jc M^*=0,
  \end{equation*}
 and 
  \begin{equation*}\label{E:4.4}
   T=\big\{\{\tz,\tf\}\in\Tmax\mid M\hz_0=0\big\}.
  \end{equation*}
\end{corollary}

If there exists $\la_0\in\Cbb$ with the property $q(\la_0)=2n$, then system~\eqref{Sla} is in the limit circle case for all 
$\la\in\Cbb$, i.e., $\qp=\qm=2n$, see Remark~\ref{R:ql.def}. Hence for any $\nu\in\Rbb$ there exist solutions 
(suppressing the argument $\nu$) $\Th^{[1]},\dots,\Th^{[2n]}$ of system \Sla{\nu}, which are linearly independent, square 
summable, and the fundamental matrix $\Th_k$ satisfies $\Th_0=I$, which implies $\Ups=\Jc$, i.e., $\rank\Ups=2n$, see 
Remark~\ref{R:Ups.def}. Upon combining the latter remark and Theorem~\ref{T:s-e.ext} we obtain the following result; cf. 
\cite[Theorem~5.10]{gR.yS14}.

\begin{corollary}\label{C:s-e.LCC}
 Let Hypothesis~\ref{H:SAC} be satisfied, assume that there exists a number $\la_0\in\Cbb$ such that $q(\la_0)=2n$, and 
 $\nu\in\Rbb$ be fixed. Let $\Th_k$ be the fundamental matrix of system~\Sla{\nu} satisfying $\Th_0=I$ and  
 denote its columns by $\Th^{[1]},\dots,\Th^{[2n]}$, i.e., $\Th_k=(\Th^{[1]}_k,\dots,\Th^{[2n]}_k)$. Then a linear relation 
 $T\subseteq \tltpt$ is  a~self-adjoint extension of $\Tmin$ if and only if there exist matrices $M,L\in\Cbb^{2n\times 2n}$ such 
 that
  \begin{equation}\label{E:4.5}
   \rank(M,L)=2n,\quad M\Jc M^*-L\,\Jc\,L^*=0,
  \end{equation}
 and 
  \begin{equation}\label{E:4.6}
   T=\Bigg\{\{\tz,\tf\}\in\Tmax\mid M \hz_0-L\msmatrix{(\Th^{[1]},\hz)_{N+1}\\ \vdots\\ (\Th^{[2n]},\hz)_{N+1}}=0\Bigg\}.
  \end{equation}
\end{corollary}

Especially, if $\Iz$ is a finite discrete interval, then the equality $q(\la)=2n$ is trivially satisfied for any $\la\in\Cbb$. 
Therefore we get from Corollary~\ref{C:s-e.LCC} yet one more special case of Theorem~\ref{T:s-e.ext}.

\begin{corollary}\label{C:s-e.finite}
 Let $\Iz$ be a finite discrete interval and Hypothesis~\ref{H:SAC} be satisfied. Then a linear relation 
 $T\subseteq \tltpt$ is a self-adjoint extension of $\Tmin$ if and only if there exist matrices 
 $M,L\in\Cbb^{2n\times 2n}$ such that
  \begin{equation}\label{E:4.7}
   \rank(M,L)=2n,\quad M\Jc\,M^*-L\,\Jc\,L^*=0,
  \end{equation}
 and 
  \begin{equation}\label{E:4.8}
   T=T_{M,L}\coloneq\big\{\{\tz,\tf\}\in\Tmax\mid M \hz_0-L\,\hz_{N+1}=0\big\}.
  \end{equation}
\end{corollary}
\begin{proof}
 By Corollary~\ref{C:s-e.LCC} every self-adjoint extension of $\Tmin$ can be expressed as in \eqref{E:4.6} with matrices
 $M,L\in\Cbb^{2n\times2n}$ satisfying \eqref{E:4.5}. If we put $\tilde{L}\coloneq L\,\Ph_{N+1}^*\,\Jc\in\Cbb^{2n\times2n}$, then 
 $M,\tilde{L}$ satisfies \eqref{E:4.7} and the linear relation in \eqref{E:4.6} can be written as $T_{M,\tilde{L}}$.
\end{proof}

One can easily observe that a linear relation $T_{M,L}$, i.e., the linear relation given by \eqref{E:4.8} with 
$M,L\in\Cbb^{2n\times2n}$ satisfying \eqref{E:4.7}, is the same as a linear relation $T_{\Mf,\Lf}$, where $\Mf\coloneq CM$ and 
$\Lf\coloneq CL$ for an arbitrary invertible matrix $C\in\Cbb^{2n\times 2n}$. We show that the converse is also true (see 
Remark~\ref{R:bound.cond.equiv}(i)). Moreover, it is well known that all self-adjoint extensions of operators associated with the 
regular second order Sturm--Liouville differential equations can be expressed by using the separated or coupled boundary 
conditions; e.g. \cite{slC.fG.rN.mZ14}. In the last part of this section we show similar results for scalar symplectic systems 
on a finite interval, i.e., $n=1$ and $N\in\Nbb$, and provide a unique representation of all self-adjoint extensions of $\Tmin$. 
The main assumptions for this treatment are summarized in the following hypothesis.

\begin{hypothesis}\label{H:5.1}
 The discrete interval $\Iz$ is finite, i.e., there exists $N\in\Nbb$ such that $\Iz=[0,N]_\sZbb$, we have $n=1$, 
 Hypothesis~\ref{H:SAC} is satisfied, and the matrices $M,L\in\Cbb^{2\times2}$ are such that \eqref{E:4.7} holds.
\end{hypothesis}

In this case, identity \eqref{E:4.7} implies either that  $\rank M=\rank L=2$, or that $\rank M=\rank L=1$, which together yield 
the following dichotomy on the boundary conditions in \eqref{E:4.8}.

\begin{theorem}\label{T:5.2}
 Let Hypothesis~\ref{H:5.1} be satisfied. Then, the following hold.
  \begin{enumerate}[leftmargin=10mm,topsep=1mm,label={\emph{(\roman*)}}]
   \item A linear relation $T_{M,L}$ given through $M,L\in\Cbb^{2\times2}$ with $\rank M=1=\rank L$ is a self-adjoint extension 
         of $\Tmin$ if and only if $T_{M,L}=T_{P,Q}\coloneq\big\{\{\tz,\tf\}\in\Tmax\mid P \hz_0=0=Q\,\hz_{N+1}\big\}$, where
          \begin{equation}\label{E:5.1}
           P=\mmatrix{\cos\al_0 & \sin\al_0\\ 0 & 0},\quad
           Q=\mmatrix{0 & 0\\ -\sin\al_{N+1} & \cos\al_{N+1}},
          \end{equation}
         for a unique pair $\al_0,\al_{N+1}\in[0,\pi)$. 
   \item A linear relation $T_{M,L}$ given through $M,L\in\Cbb^{2\times2}$ with $\rank M=2=\rank L$ is a self-adjoint extension 
         of $\Tmin$ if and only if $T_{M,L}=T_{R,\be}\coloneq\big\{\{\tz,\tf\}\in\Tmax\mid \e^{i\be}R\,\hz_0=\hz_{N+1}\big\}$ 
         with a unique $\be\in[0,\pi)$ and a symplectic matrix $R\in\Rbb^{2\times2}$.
  \end{enumerate}
\end{theorem}
\begin{proof}
 Since the pairs of matrices $P,Q$ and $\e^{i\be}R, I$ satisfy \eqref{E:4.7}, Corollary~\ref{C:s-e.finite} implies that the 
 linear relations $T_{P,Q}$ and $T_{R,\be}$ are self-adjoint extensions of $\Tmin$.
 
 (i) Let $T_{M,L}$ be a linear relation given through $M,L\in\Cbb^{2\times2}$ satisfying \eqref{E:4.7} and with 
    $\rank M=1=\rank L$. Since by \eqref{E:2.3} we have $\dim[\Rc(M)\cap\Rc(L)]=0$, it follows that $M\xi=L\eta$ for some 
    $\xi,\eta\in\Cbb^2$ if and only if $M\xi=0=L\eta$. Therefore, the boundary conditions in \eqref{E:4.8} can be expressed as 
    $M\hz_0=0=L\hz_{N+1}$. The rank condition implies that $M=ab^\top$ and $L=cd^\top$ for some vectors 
    $a,b,c,d\in\Cbb^2\stm\{0\}$. Then the equality $M\Jc M^*=0=L\Jc L^*$ does not depend on the vectors $a,c$ and it is 
    equivalent with $b^\top\!\Jc\,b=0=d^\top\!\Jc d$, which implies that $b$ and $d$ are (scalar) complex multiples of vectors 
    from $\Rbb^2$. Therefore, without loss of generality, $a,c$ may be chosen such that $M,L$ can be written in the form as in 
    \eqref{E:5.1} for some $\al_0,\al_{N+1}\in[0,\pi)$. The uniqueness follows from the fact that $\cotan\al=\cotan\be$ with 
    $\al,\be\in(0,\pi)$ if and only if $\al=\be$.
 
 (ii) Finally, let $T_{M,L}$ be a linear relation given through $M,L\in\Cbb^{2\times2}$ satisfying \eqref{E:4.7} and with 
      $\rank M=2=\rank L$. Then the boundary conditions in \eqref{E:4.8} can be written as $\hz_{N+1}=K\hz_0$, where 
      $K\coloneq L^{-1}M$. Upon applying the second equality in \eqref{E:4.7} we obtain that the matrix $K$ is conjugate 
      symplectic, i.e., $K\Jc K^*=\Jc$. Therefore, $K^{-1}=-\Jc K^*\Jc$ and $\abs{\det K}=1$, i.e., $\det K=\e^{i\de}$ for some 
      $\de\in[0,2\pi)$, which implies $K^{-1}=\e^{-i\de}\adj(K)=-\e^{i\de}\Jc K^\top\!\Jc$, i.e., 
      $K^{*\top}=\overline{K}=\e^{i\de} K$. If we put $R\coloneq \e^{-i\de/2} K$, i.e., $K=\e^{i\de/2} R$, then $\overline{R}=R$ 
      and $\det R=1$, i.e., $R\in\Rbb^{2\times2}$ is a symplectic matrix. Uniqueness can be verified by a~direct calculation.
\end{proof}

As an illustration of the last theorem we provide a description of the Krein--von~Neumann extension of the minimal linear 
relation $\Tmin$ under Hypothesis~\ref{H:5.1}.

\begin{example}\label{Ex:krein.ext}
 Assume that system~\eqref{Sla} is such that Hypothesis~\ref{H:5.1} holds and that the minimal linear relation $\Tmin$ is 
 positive, i.e., there exists $c>0$ such that $\innerP{\tz}{\tf}\geq c\,\normP{\tz}$ for all $\{\tz,\tf\}\in\Tmin$. Then the 
 Krein--von~Neumann self-adjoint extension extension of $\Tmin$ admits the representation given in \eqref{E:krein.ext}, i.e.,
  \begin{equation*}
   T_K=\Tmin\dotplus(\ker\Tmax\times\{0\}).
  \end{equation*}
 We show that $T_K$ can be also expressed as in the second part of Theorem~\ref{T:5.2} with a suitable matrix $R$ and a~number 
 $\be\in[0,2\pi)$. By definition,
  \begin{equation*}
   \ker\Tmax=\{\tz\in\ltp\mid \{\tz,\tilde{0}\}\in\Tmax\},
  \end{equation*}
 i.e., $\hz$ solves \Sla{0}, i.e., $\mL(\hz)_k=0$ on $[0,N]_\sZbb$. Because all solutions of \Sla{0} are square summable in 
 this case, Hypothesis~\ref{H:SAC} implies that $\dim\ker\Tmax=2$. If $\tz\in\dom T_K$, then there exist $\ty\in\dom\Tmin$ and 
 $\tw\in\ker\Tmax$ such that $\tz=\ty+\tw$ or
  \begin{equation}\label{E:5.4a}
   \hz_k=\hy_k+\hw_k\quad \text{for all $k\in[0,N+1]_\sZbb$},
  \end{equation}
 where $\hz\in\tz$, $\hy\in\ty$, and $\hw\in\tw$ are the uniquely determined elements. Moreover, $\hy_0=0=\hy_{N+1}$ by 
 \eqref{E:2.40} and $\hw_k=\al^{[1]}\hw^{[1]}_k+\al^{[2]}\hw^{[2]}_k$ for all $k\in[0,N+1]_\sZbb$, where $\hw^{[1]}$ and 
 $\hw^{[2]}$ form a basis of $\ker\Tmax$.
 
 Let us define the matrix 
 $\Gc=\msmatrix{a & b\\ c & d}\coloneq (\Sc_0\times\Sc_1\times \dots\times \Sc_N)^{-1}\in\Cbb^{2\times2}$. Then one easily  
 concludes that the matrix $\Gc$ is symplectic and every solution $z\in\Cbb([0,N+1]_\sZbb)^{2}$ of system \Sla{0} satisfies
  \begin{equation}\label{E:5.4b}
   z_{N+1}=\Gc z_0.
  \end{equation}
 In the following construction we consider two cases: either $b\neq0$ or $b=0$. 
 
First, assume that 
 $b\neq0$. Then there exist two solutions of system \Sla{0} such that
  \begin{equation*}
   \hw^{[1]}_0=\mmatrix{0\\ 1/b},\quad 
   \hw^{[2]}_0=\mmatrix{1\\ -a/b}.
  \end{equation*}
 These solutions are obviously linearly independent and by \eqref{E:5.4b} we have
  \begin{equation*}
   \hw^{[1]}_{N+1}=\mmatrix{1\\ a/b},\quad 
   \hw^{[2]}_{N+1}=\mmatrix{0\\ c-da/b}.
  \end{equation*}
 If we take these two solutions as a basis of $\ker\Tmax$, then \eqref{E:5.4a} yields
  \begin{equation*}
   \hz_k=\hy_k+\al^{[1]}\hw^{[1]}_k+\al^{[2]}\hw^{[2]}_k\quad \text{for all $k\in[0,N+1]_\sZbb$}.
  \end{equation*}
 Upon evaluating $\hz_k$ at $k=0$ and $k=N+1$ we obtain
  \begin{equation*}
   \hz_0=\mmatrix{\al^{[2]}\\ \al^{[1]}/b-\al^{[2]}a/b},\quad
   \hz_{N+1}=\mmatrix{\al^{[1]}\\ \al^{[1]}d/b+\al^{[2]}c-\al^{[2]}da/b},
  \end{equation*}
 which for $\hz_k=\msmatrix{\hx_k\\ \hu_k}$ implies $\al^{[1]}=\hx_{N+1}$ and $\al^{[2]}=\hx_0$. Therefore,
  \begin{equation*}
   \mmatrix{\hx_{N+1}\\ \hx_{N+1}\,d/b+\hx_0\,c-\hx_0\,da/b}=\hz_{N+1}
   =\Gc\,\hz_0=\Gc\mmatrix{\hx_0\\ \hx_{N+1}/b-\hx_0\,a/b}.
  \end{equation*}
 It means that $\hz\in\dom T_{R,\be}$, where $\be\in[0,\pi)$ is such that $\e^{i\be}=\sqrt{ad-bc}$, and $R=\e^{-i\be}\,\Gc$, 
 i.e., $T_K\subseteq T_{R,\be}$. On the other hand, $T_K$ and $T_{R,\be}$ are self-adjoint extensions of $\Tmin$, thus 
 $T_K=T_{R,\be}$. Especially, if the coefficients $a,b,c,d$ are real, then $T_{R,\be}=T_{\Gc,0}$.
 
 If $b=0$, then $\Gc=\msmatrix{a & 0\\ c & d}$ with $\abs{ad}=1$, i.e., $d\neq0$. In this case we proceed in the same way 
 with the basis of $\ker\Tmax$ given by the solutions $\tw^{[1]}$ and $\tw^{[2]}$ of \Sla{0} such that
  \begin{equation*}
   \hw^{[1]}_0=\mmatrix{0\\ 1/d},\quad 
   \hw^{[2]}_0=\mmatrix{1\\ -c/d}.
  \end{equation*}
 Then $\msmatrix{\hx_0\,a\\ \hu_{N+1}}=\hz_{N+1}=\Gc\,\hz_0=\Gc\msmatrix{\hx_0\\ \hu_{N+1}/d - \hx_0\,c/d}$. This shows 
 (again) that $T_K=T_{R,\be}$ with $\be\in[0,\pi)$ being such that $\e^{i\be}=\sqrt{ad}$, and $R=\e^{-i\be}\,\Gc$.

 In particular, let $\Sc_k=\msmatrix{1 & -b_k\\ 0 & 1}$ and $\Ps_k=\msmatrix{w_k & 0\\ 0 & 0}$ with $b_k>0$ and 
 $w_k>0$ on  $[0,N]_\sZbb$. This system satisfies Hypothesis~\ref{H:SAC} and corresponds to the second order Sturm--Liouville 
 difference equation $-\De[p_k\,\De y_{k-1}(\la)]=\la\,w_k\,y_k(\la)$ with $b_k=1/p_{k+1}$ (see Example~\ref{Ex:2n.order}(i)). 
 Then $\Gc=\msmatrix{1 & \sum_{k=0}^N b_k\\ 0 & 1}$ and by the previous part we have
  \begin{equation*}
   T_K=\Big\{\{\tz,\tf\}\in\Tmax\mid \hz=\msmatrix{\hx\\ \hu}\in\Cbb([0,N+1]_\sZbb)^{2},\ 
                                \hu_0=\hu_{N+1}=\Big(\sum_{k=0}^N b_k\Big)^{\!-1}\times(\hx_{N+1}-\hx_0)\Big\}.\qedhere
  \end{equation*}
\end{example}

The boundary conditions in Theorem~\ref{T:5.2} include four particular cases. Namely, for $\al_0=0$ 
and $\al_{N+1}=\pi/2$ we get the Dirichlet boundary conditions $\hx_0=0=\hx_{N+1}$, while for $\al_0=\pi/2$ 
and $\al_{N+1}=0$ we have the Neumann boundary conditions $\hu_0=0=\hu_{N+1}$, where $\hz_k=\msmatrix{\hx_k\\ \hu_k}$. The choice 
$R=I$ and $\be=0$ yields the periodic boundary conditions $\hz_0=\hz_{N+1}$ and the choice $R=I$ and $\be=\pi$ leads to the 
antiperiodic boundary conditions $\hz_0=-\hz_{N+1}$.

In the first part of the following theorem we show that any self-adjoint extension of $\Tmin$ can be described by using the 
matrices determining the Dirichlet and Neumann boundary conditions. For convenience, we introduce the general boundary trace map 
$\ga_{M,L}:\Cbb(\Izp)^{2}\to\Cbb^{2}$ as
 \begin{equation*}
  \ga_{M,L}(\hz)\coloneq M\hz_0-L\hz_{N+1}, 
 \end{equation*}
see also \cite{slC.fG.rN.mZ14}. Then $T_{M,L}=\big\{\{\tz,\tf\}\in\Tmax\mid \ga_{M,L}(\hz)=0\big\}$. Especially, for $P,Q$ 
given in \eqref{E:5.1} we denote $\ga_x\coloneq\ga_{P,Q}$ for $\al_0=0$, $\al_{N+1}=\pi/2$, i.e., $\ga_x(\hz)=0$ abbreviates the 
Dirichlet boundary conditions, and similarly $\ga_u\coloneq\ga_{P,Q}$ for $\al_0=\pi/2$, $\al_{N+1}=0$, i.e., $\ga_u(\hz)=0$ 
abbreviates the Neumann boundary conditions. In the second part of this theorem we derive yet another equivalent representation 
of $T_{M,L}$, which possesses the uniqueness property.

\begin{theorem}\label{C:5.3}
 Let Hypothesis~\ref{H:5.1} be satisfied. Then the following hold.
  \begin{enumerate}[leftmargin=10mm,topsep=1mm,label={\emph{(\roman*)}}]
   \item A linear relation $T$ is a self-adjoint extension of $\Tmin$ if and only if there exist matrices $F,G\in\Cbb^{2\times2}$ 
         such that
          \begin{equation}\label{E:5.10}
           \rank(F,G)=2,\quad FG^*=GF^*
          \end{equation}
         and 
          \begin{equation}\label{E:5.11}
           T=T_{F,G}\coloneq\big\{\{\tz,\tf\}\in\Tmax\mid F\,\ga_x(\hz)+G\,\ga_u(\hz)=0\big\}.
          \end{equation}
   \item We have $T_{F,G}=T_{\Ff,\Gf}$, where $\Ff,\Gf$ satisfy \eqref{E:5.10}, if and only if $\Ff=CF$ and $\Gf=CG$ for some 
         invertible matrix $C\in\Cbb^{2\times2}$.
   \item A linear relation $T$ is a self-adjoint extension of $\Tmin$ if and only if there exists a unitary matrix 
         $V\in\Cbb^{2\times2}$ such that
          \begin{equation}\label{E:5.13}
           T=T_{V}\coloneq\big\{\{\tz,\tf\}\in\Tmax\mid i(V-I)\,\ga_x(\hz)=(V+I)\,\ga_u(\hz)\big\}.
          \end{equation}
   \item We have $T_{V}=T_{\Vf}$, where $\Vf\in\Cbb^{2\times2}$ is a unitary matrix, if and only if $\Vf=V$.
  \end{enumerate}
\end{theorem}
\begin{proof}
 (i) Let $T$ be given by \eqref{E:5.11} with $F,G\in\Cbb^{2\times2}$ satisfying \eqref{E:5.10}. If we put 
     $M\coloneq FP_0+GP_{\pi/2}$ and $L\coloneq FQ_{\pi/2}+GQ_{0}$, where $P_{\om}$ and $Q_{\om}$ are the matrices corresponding 
     to $P,Q$ defined in \eqref{E:5.1} with $\om\in\{0,\pi/2\}$. Then $M\Jc M^*-L\,\Jc L^*=FG^*-GF^*=0$ and $\rank(F,G)=2$ is 
     equivalent with $\rank(M,L)=2$. Hence $M,L$ satisfy \eqref{E:4.7}. Moreover, for the left-hand side of the boundary 
     conditions in \eqref{E:5.11} we have $F\,\ga_x(\hz)+G\,\ga_u(\hz)=\ga_{M,L}(\hz)$. Therefore $\{\tz,\tf\}\in T_{M,L}$ if 
     and only if $\{\tz,\tf\}\in T_{F,G}$, i.e., $T_{F,G}$ is a self-adjoint extension of $\Tmin$ by Corollary~\ref{C:s-e.finite}.
     On the other hand, let $T$ be a self-adjoint extension of $\Tmin$, i.e., $T=T_{M,L}$ with $M,L\in\Cbb^{2\times2}$ 
     satisfying \eqref{E:4.7}. If we put $F\coloneq MP_0-LP_{\pi/2}$ and $G\coloneq LQ_0-MQ_{\pi/2}$, then the 
     conditions in \eqref{E:5.10} hold and $\ga_{M,L}(\hz)$ can be written as in \eqref{E:5.11}.
     
 (ii) Sufficiency is clear. Assume that $T_{F,G}=T_{\Ff,\Gf}$ for two pairs of  matrices $F,G$ and $\Ff,\Gf$ satisfying 
      \eqref{E:5.10}. Then, by \eqref{E:5.11}, we have for any $\{\tz,\tf\}\in\Tmax$ that $F\,\ga_x(\hz)+G\,\ga_u(\hz)=0$ if 
      and only if $\Ff\,\ga_x(\hz)+\Gf\,\ga_u(\hz)=0$. It means that $\hz_0,\hz_{N+1}$ solve simultaneously the both systems of 
      algebraic equations with the coefficient matrices $F,G$ and $\Ff,\Gf$. It means that these systems are equivalent, which 
      implies an existence of an invertible matrix $C\in\Cbb^{2\times2}$ such that $\Ff=CF$ and $\Gf=CG$.
      
 (iii) Let $T$ be given by \eqref{E:5.13} with a unitary matrix $V\in\Cbb^{2\times2}$. If we put $F\coloneq\frac{i}{2}(I-V)$ and 
       $G\coloneq \frac{1}{2}(I+V)$. Then $FG^*=GF^*$ and, by \eqref{E:2.2}, $\rank(F,G)=2$, i.e., $F,G$ satisfy \eqref{E:5.10}. 
       Since the boundary conditions in \eqref{E:5.13} are equivalent with the boundary conditions in \eqref{E:5.11} with $F,G$ 
       defined above, i.e., $\{\tz,\tf\}\in T_{F,G}$ if and only if $\{\tz,\tf\}\in T_{V}$, it follows from the previous part 
that 
       the linear relation $T_V$ is a self-adjoint extension of $\Tmin$. On the other hand, let $T$ be a self-adjoint extension 
       of $\Tmin$. Then, by the part (i), we have $T=T_{F,G}$ with $F,G\in\Cbb^{2\times2}$ satisfying \eqref{E:5.10}. Since by 
       \eqref{E:2.2} and \eqref{E:5.10} we have $\rank(F+iG)=2$, the matrix $V\coloneq (F+iG)^{-1}(iG-F)$ is well-defined. One 
       can directly verify that $V$ is a unitary matrix and the boundary conditions $F\,\ga_x(\hz)+G\,\ga_u(\hz)=0$ are satisfied 
       if and only if $i(V-I)\,\ga_x(\hz)-(V+I)\,\ga_u(\hz)=0$, i.e., $T_{F,G}=T_V$.
       
 (iv) If $V=\Vf$, then $T_V=T_\Vf$. On the other hand, assume that $T_V=T_\Vf$ for two unitary matrices 
      $V,\Vf\in\Cbb^{2\times2}$. Then $T_{F,G}=T_V=T_\Vf=T_{\Ff,\Gf}$ with $F,G$ and $\Ff,\Gf$ being given as in the previous 
      part. Then $V=(F+iG)^{-1}(iG-F)$ and $\Vf=(\Ff+i\Gf)^{-1}(i\Gf-\Ff)$ and by the part (ii) there exists an invertible matrix 
      $C\in\Cbb^{2\times2}$ such that $\Ff=CF$ and $\Gf=CG$. Upon combining these facts we obtain $V=\Vf$.  
\end{proof}

\begin{remark}\label{R:bound.cond.equiv} \
 \begin{enumerate}[leftmargin=10mm,topsep=1mm]
  \item As a consequence of Theorem~\ref{T:5.2}(i)-(ii) we obtain that $T_{M,L}=T_{\Mf,\Lf}$ if and only if $\Mf=CM$ and $\Lf=CL$ 
        for some invertible matrix $C\in\Cbb^{2\times2}$.
  \item The statement of Theorem~\ref{T:5.2}(iii)-(iv) shows that the map from the set of all $2\times2$ unitary matrices 
        to the set of all self-adjoint extensions expressed as in \eqref{E:5.13} is a bijection.
 \end{enumerate}
\end{remark}

\section{Proof of main result}\label{S:proof}

In this section, a proof is given for Theorem~\ref{T:s-e.ext} which utilizes several arguments from the linear algebra and whose 
main idea goes back to \cite{jS86}. It is based on a construction of a suitable GKN-set (see Theorem~\ref{T:general.s-e.ext}), 
and on a more convenient expression than that given in \eqref{E:3.5}  for elements in $\dom\Tmax$. Similar results for 
system~\eqref{E:1.2.hamilton} can be found in \cite[Section~4]{gR.yS14}. 

\begin{lemma}\label{L:3.2}
 Let Hypothesis~\ref{H:SAC} be satisfied, $\{\tz,\tf\}\in\Tmax$ arbitrary, and $\vp^{[1]},\dots,\vp^{[\qp]}$ be arranged such 
 that \eqref{E:rank.Om.1.2} holds. Then the element $\hz$ can be uniquely expressed as
  \begin{equation}\label{E:3.8}
   \hz_k=\hy_k+\sum_{i=1}^{2n} \eta_i\,\hz_k^{[i]}+\sum_{j=1}^{p-2n} \zeta_j\,\vp_k^{[j]},\quad k\in\Izp,
  \end{equation}
 where $\hy\in\dom\Tmin$, $\hz^{[1]},\dots,\hz^{[2n]}$ are specified in Lemma~\ref{L:2.6}, and $\eta_i,\zeta_j\in\Cbb$ 
 for all $i\in\{1,\dots,2n\}$ and $j\in\{1,\dots,p-2n\}$. Moreover,
  \begin{equation}\label{E:3.7}
   \rank \Om_{p-2n}=p-2n,
  \end{equation}
 where $\Om$ was defined in \eqref{E:Om.def}.
\end{lemma}
\begin{proof}
 Since \eqref{E:rank.Om.1.2} is satisfied, there exists an invertible matrix $P\in\Cbb^{p\times p}$ such that
  \begin{equation}\label{E:3.16}
   \Om P=\mmatrix{I_{p-2n} & 0\\ Q & R},
  \end{equation}
 where $I_{p-2n}$ is the $(p-2n)\times(p-2n)$ identity matrix and $0$ stands for the $(p-2n)\times2n$ zero matrix. If we put 
 $\Xi=(\Xi^{[1]},\Xi^{[2]})$, where $\Xi^{[1]}\in\Cbb^{2n\times(p-2n)}$ and $\Xi^{[2]}\in\Cbb^{2n\times2n}$, and multiply 
 \eqref{E:3.16} by $\overline{\Xi}$ from the left, we obtain
  \begin{equation*}\label{E:3.17}
   \overline{\Xi^{[1]}}=-\overline{\Xi^{[2]}}\,Q,
  \end{equation*}
 i.e., $\overline{\Xi}=\big(-\overline{\Xi^{[2]}}\,Q,\overline{\Xi^{[2]}}\big)$. It implies that $\rank\Xi^{[2]}=2n$ by the 
 second inequality in \eqref{E:2.1}, because $\rank\Xi=2n$. By multiplying equality \eqref{E:3.6A} by the matrix 
 $(\Xi^{[2]})^{\top-1}$ from the right, we get
  \begin{equation*}\label{E:3.18}
   Z_k\,(\Xi^{[2]})^{\top-1}=Y_k\,(\Xi^{[2]})^{\top-1}+\Ph_k^{[1]}\,\Xi^{[1]\top}(\Xi^{[2]})^{\top-1}+\Ph_k^{[2]},
  \end{equation*}
 where $\Ph_k^{[1]}\in\Cbb^{2n\times(p-2n)}$ and $\Ph_k^{[2]}\in\Cbb^{2n\times2n}$ are such that 
 $\Ph_k=(\Ph_k^{[1]},\Ph_k^{[2]})$.
 It shows that every solution $\vp^{[2n-p+1]},\dots,\vp^{[p]}$ can be uniquely expressed with $\hy^{[i]}$, $\hz^{[i]}$, 
 $i\in\{1,\dots,2n\}$, and $\vp^{[1]},\dots,\vp^{[p-2n]}$, i.e., 
  \begin{equation}\label{E:3.20}
   \vp_k^{[j]}=\hu_k^{[j]}+\sum_{r=1}^{2n}\eta_{j,r}\,\hz_k^{[r]}+\sum_{s=1}^{p-2n}\zeta_{j,s}\,\vp_k^{[s]},
   \quad k\in\Izp,\quad j\in\{p-2n+1,\dots,p\},
  \end{equation}
 for some $\hu_k^{[j]}\in\dom\Tmin$ and $\eta_{j,r},\zeta_{j,s}\in\Cbb$. Therefore, the expression in~\eqref{E:3.8} follows from 
 \eqref{E:3.5}. Moreover, if we multiply \eqref{E:3.20} by $\vp_k^{[i]*}\Jc$ from the left, where $i\in\{1,\dots,p-2n\}$, then 
  \begin{equation*}\label{E:3.22}
   (\vp^{[i]},\vp^{[j]})_{N+1}=(\vp^{[i]},\hy^{[j]})_{N+1}+\sum_{r=1}^{2n}\eta_{j,r}\,(\vp^{[i]},\hz^{[r]})_{N+1}
    +\sum_{s=1}^{p-2n}\zeta_{j,s}\,(\vp^{[i]},\vp^{[s]})_{N+1}.
  \end{equation*}
 Hence from \eqref{E:2.39} and the definition of $\hz^{[i]}$ we have
  \begin{equation}\label{E:3.23}
   \Om^{[1,2]}_{p-2n,\,\qm}=\Om_{p-2n}\,T^\top,
  \end{equation}
 where $T\in\Cbb^{\qm\times(p-2n)}$ is a matrix consisting of the elements $\zeta_{j,s}$ for $j\in\{\qp+1,\dots,p\}$ and 
 $s\in\{1,\dots,p-2n\}$. Since the solutions are arranged such that $\rank\Om_{p-2n,\,\qm}^{[1,2]}=p-2n$, identity \eqref{E:3.7} 
 follows from \eqref{E:3.23} and the second inequality in \eqref{E:2.1}.
\end{proof}

\begin{remark}\label{R:3.3}
 If we switch the role of $v^{[\cdot]}(\la_0)$ and $w^{[\cdot]}(\bla_0)$ in the definition of $\vp^{[1]},\dots,\vp^{[p]}$ in 
 \eqref{E:3.3+3.4}, i.e., we put $\vp^{[i]}=w^{[i]}(\bla_0)$ for $i\in\{1,\dots,\qm\}$ and 
 $\vp^{[j+\qm]}=v^{[j]}(\la_0)$ for $j\in\{1,\dots,\qp\}$, then the solutions $\vp^{[1]},\dots,\vp^{[\qm]}$ can be arranged such 
 that \eqref{E:3.8} and \eqref{E:3.7} hold.
\end{remark}

Now, we give the proof of Theorem~\ref{T:s-e.ext}.

\begin{proof}[Proof of Theorem~\ref{T:s-e.ext}]
 Assume that $T$ is a self-adjoint extension of $\Tmin$. Then, by Theorem~\ref{T:general.s-e.ext} there exists a GKN-set 
 $\{\be_j\}_{j=1}^q$ for $(\Tmin,\Tmax)$ such that \eqref{E:2.32} holds. Since $\be_j\in\Tmax$, they may be identified as
 $\be_j=\{\tw^{[j]},\th^{[j]}\}\in\Tmax$. By Lemma~\ref{L:3.2}, the elements $\hw^{[j]}$ can be uniquely expressed as
  \begin{equation}\label{E:4.9}
   \hw_k^{[j]}=\hy_k^{[j]}+\sum_{i=1}^{2n}\eta_{j,i}\,\hz_k^{[i]}+\sum_{l=1}^{2q-2n}\zeta_{j,l}\,\vp_k^{[l]},\quad k\in\Izp,
  \end{equation}
 where $\hy^{[j]}\in\dom\Tmin$ and $\eta_{j,i},\zeta_{j,l}\in\Cbb$. We  next show that the matrices
  \begin{equation*}\label{E:4.10}
   M\coloneq (\hw_0^{[1]},\dots,\hw_0^{[q]})^*\Jc\in\Cbb^{q\times2n},\quad
   L\coloneq \mmatrix{\overline{\zeta_{1,1}} & \cdots & \overline{\zeta_{1,2q-2n}}\\
                      \vdots & \ddots & \vdots\\
                      \overline{\zeta_{q,1}} & \cdots & \overline{\zeta_{q,2q-2n}}}\in\Cbb^{q\times(2q-2n)}
  \end{equation*}
 satisfy \eqref{E:4.1}. 
 
 Since $\rank(M,L)\leq q$, assume that $\rank(M,L)<q$. Then, there exists 
 $C=(c_1,\dots,c_q)^\top\in\Cbb^{q}\stm\{0\}$ such that $C^*(M,L)=0$, i.e., $C^*M=0=C^*L$. If 
 $\hw_k\coloneq\sum_{j=1}^q c_j \hw_k^{[j]}$ for $k\in\Izp$, then $\hw_0=\Jc M^*C=0$ and also 
 $(\hw,\vp^{[i]})_{N+1}=\sum_{j=1}^q \overline{c_j}\,(\hw^{[j]},\vp^{[i]})_{N+1}$ for all $i\in\{1,\dots,2q-2n\}$. Hence by 
 \eqref{E:4.9} and \eqref{E:2.39} we have
  \begin{equation*}\label{E:4.13}
   \big((\hw,\vp^{[1]})_{N+1},\dots,(\hw,\vp^{[2q-2n]})_{N+1}\big)=C^*L\,\Om_{2q-2n}=0.
  \end{equation*}
 But then $(\hw,\hy)_{N+1}=0$ for any $\hy\in\dom\Tmax$, because it can be written as in \eqref{E:3.8}. It means that 
 $\hw\in\dom\Tmin$ by \eqref{E:2.39} and hence $\be_1,\dots,\be_q$ are linearly dependent in $\Tmax$ modulo $\Tmin$, 
 which contradicts the assumption that that $\{\be_j\}_{j=1}^q$ is a GKN-set. Therefore, the first condition in \eqref{E:4.1} is 
 satisfied. 
 
 Next, we see that
  \begin{equation}\label{E:4.14}
   \mmatrix{(\hw^{[1]},\hw^{[1]})_0 & \cdots & (\hw^{[1]},\hw^{[q]})_0\\
             \vdots & \ddots & \vdots\\
             (\hw^{[q]},\hw^{[1]})_0 & \hdots & (\hw^{[q]},\hw^{[q]})_0}=M\Jc M^*  
  \end{equation}
 and by using \eqref{E:4.9}, \eqref{E:2.39}, and the definition of $\hz^{[i]}$, also see that
  \begin{equation}\label{E:4.15}
   \mmatrix{(\hw^{[1]},\hw^{[1]})_{N+1} & \cdots & (\hw^{[1]},\hw^{[q]})_{N+1}\\
             \vdots & \ddots & \vdots\\
             (\hw^{[q]},\hw^{[1]})_{N+1} & \hdots & (\hw^{[q]},\hw^{[q]})_{N+1}}=L\,\Om_{2q-2n}\,L^*.
  \end{equation}
 Since $\{\be_j\}_{j=1}^q$ is a GKN-set, we obtain from \eqref{E:2.38} that 
  \begin{equation*}\label{E:4.16}
   0=[\be_i:\be_j]=(\hw^{[i]},\hw^{[j]})_k\big|_{0}^{N+1}
  \end{equation*}
 for all $i,j\in\{1,\dots,q\}$. By \eqref{E:4.14} and \eqref{E:4.15}, this implies that $M\Jc M^*-L\,\Om_{2q-2n}\,L^*=0$,  and 
 that the second condition in \eqref{E:4.1} is also satisfied. 
 
 For any $\hz\in\dom\Tmax$, we can write
  \begin{equation}\label{E:4.17+18}
   \mmatrix{(\hw^{[1]},\hz)_0\\ \vdots\\ (\hw^{[q]},\hz)_0}=M\hz_0,\qquad
   \mmatrix{(\hw^{[1]},\hz)_{N+1}\\ \vdots\\ (\hw^{[q]},\hz)_{N+1}}
     =L\!\mmatrix{(\vp^{[1]},\hz)_{N+1}\\ \vdots\\ (\vp^{[2q-2n]},\hz)_{N+1}},
  \end{equation}
 where the second equality follows from \eqref{E:4.9}, \eqref{E:2.39}, and the definition of $\hz^{[i]}$. Upon combining
 \eqref{E:2.32}, \eqref{E:2.38}, \eqref{E:4.17+18}, we obtain that $T$ can be expressed as
  \begin{align*}
   T &=\big\{\{\tz,\tf\}\in\Tmax\mid (\hz,\hw^{[j]})_k\big|_0^{N+1}=0\ \text{for all } j=1,\dots,q\big\}\\
     &=\big\{\{\tz,\tf\}\in\Tmax\mid \hw^{[j]*}_k\Jc\hz_k\big|_0^{N+1}=0\ \text{for all } j=1,\dots,q\big\}\\
     &=\Bigg\{\{\tz,\tf\}\in\Tmax\mid M \hz_0-L\msmatrix{(\vp^{[1]},\hz)_{N+1}\\ \vdots\\ (\vp^{[2q-2n]},\hz)_{N+1}}=0\Bigg\},
  \end{align*}
 i.e., as written in \eqref{E:4.2}.
 
 On the other hand, let $M\in\Cbb^{q\times2n}$ and $L\in\Cbb^{q\times(2q-2n)}$ satisfy \eqref{E:4.1} and $T$ be given by 
 \eqref{E:4.2}. We then must show that there exists a GKN-set $\{\be_j\}_{j=1}^q$ for $(\Tmin,\Tmax)$ such that $T$ can be 
 expressed as in~\eqref{E:2.32}. Denote the columns of $\Jc M^*\in\Cbb^{2n\times q}$ as $\rho_{1},\dots,\rho_{q}$ and the columns 
 of the matrix $(\vp^{[1]}_k,\dots,\vp^{[2q-2n]}_k)\,L^*\in\Cbb^{2n\times q}$ as $w_k^{[1]},\dots, w_k^{[q]}$, i.e.,
  \begin{equation}\label{E:4.21}
   \rho_{i}\coloneq \Jc M^* e_i,\quad 
   w_k^{[i]}\coloneq \sum_{l=1}^{2q-2n}\overline{\eta_{i,l}}\,\vp^{[l]}_k,\quad i\in\{1,\dots,q\},
  \end{equation}
 where $e_i$ is the $i$-th canonical unit vector in $\Cbb^q$ and $\eta_{i,j}$ are the elements of $L$ for $i\in\{1,\dots,q\}$ 
 and $j\in\{1,\dots,2q-2n\}$. Then, $w^{[i]}\in\Tmax$ for all $i\in\{1,\dots,q\}$ and, by Lemma~\ref{L:2.6}, there exist 
 $\be_i\coloneq\{\ty^{[i]},\th^{[i]}\}\in\Tmax$ such that
  \begin{equation*}\label{E:4.22}
   \hy^{[i]}_0=\rho_{i},\quad \hy^{[i]}_k=w^{[i]}_k,\quad k\in[b+1,\infty)_\sZbb\cap\Izp
  \end{equation*}
 for all $i\in\{1,\dots,q\}$, where the number $b$ is determined in Hypothesis~\ref{H:SAC}. We next show that $\{\be_i\}_{i=1}^q$ 
 form a GKN-set for $(\Tmin,\Tmax)$. 
 
 Since the linear independence of $\be_1,\dots,\be_q$ in $\Tmax$ modulo $\Tmin$ is equivalent to the linear independence of 
 $\hy^{[1]},\dots,\hy^{[q]}$ in $\dom\Tmax$ modulo $\Tmin$, we assume that there exists 
 $C=(c_1,\dots,c_q)^\top\in\Cbb^q\stm\{0\}$ such that
  \begin{equation*}\label{E:4.23}
   \hy\coloneq\sum_{j=1}^q c_j\,\hy^{[j]}\in\dom\Tmin.
  \end{equation*}
 Then, from \eqref{E:2.39} and \eqref{E:4.21}, we have for all $\vp^{[1]},\dots,\vp^{[2q-2n]}\in\Tmax$ that
  \begin{equation*}\label{E:4.24}
   0=\big((\hy,\vp^{[1]})_{N+1},\dots,(\hy,\vp^{[2q-2n]})_{N+1}\big)=C^* L\,\Om_{2q-2n}.
  \end{equation*}
 This implies $C^*L=0$, because $\Om_{2q-2n}$ is assumed to be invertible. Simultaneously we have $\hy_0=0$, which yields
  \begin{equation*}\label{E:4.25}
   0=\hy_0=\sum_{j=1}^q c_j\,\hy^{[j]}_0=\Jc M^* C,
  \end{equation*}
 i.e., $C^* M=0$, because the matrix $\Jc$ is invertible. But this means $C^*(M,L)=0$, which contradicts the first assumption in 
 \eqref{E:4.1}. 
 
 Next, let
  \begin{equation*}
   Y_k\coloneq \mmatrix{(\hy^{[1]},\hy^{[1]})_k & \cdots & (\hy^{[1]},\hy^{[q]})_k\\
                       \vdots & \ddots & \vdots\\
                       (\hy^{[q]},\hy^{[1]})_k & \cdots & (\hy^{[q]},\hy^{[q]})_k}.
  \end{equation*}
 Since it can be directly calculated that $Y_0=M\Jc M^*$ and $Y_{N+1}=L\,\Om_{2q-2n} L^*$, the second equality in 
 \eqref{E:4.1} implies $Y_0-Y_{N+1}=0$. Therefore, by using \eqref{E:2.38}, we get
  \begin{equation*}
   [\be_i:\be_j]=(\hy^{[i]},\hy^{[j]})_k\big|_{0}^{N+1}=0,
  \end{equation*}
 which shows that $\{\be_i\}_{i=1}^q$ is a GKN-set for $(\Tmin,\Tmax)$ as defined in Subsection~\ref{Ss:relations}. 
 
 Finally, let $\{\tw,\tg\}\in\Tmax$ be arbitrary, then
  \begin{equation}\label{E:4.29+30}
   M\hw_0=\mmatrix{(\hy^{[1]},\hw)_{0}\\ \vdots\\ (\hy^{[q]},\hw)_{0}},\quad
   L\!\mmatrix{(\vp^{[1]},\hw)_{N+1}\\ \vdots\\ (\vp^{[2q-2n]},\hw)_{N+1}}=
      \mmatrix{(\hy^{[1]},\hw)_{N+1}\\ \vdots\\ (\hy^{[q]},\hw)_{N+1}}.
  \end{equation}
 By \eqref{E:2.38} the condition $[\{\tw,\tg\}:\be_i]=0$ is equivalent to
  \begin{equation}\label{E:4.30a}
   (\hw,\hy^{[i]})_k\big|_0^{N+1}=0=-(\hy^{[i]},\hw)_k\big|_0^{N+1}
  \end{equation}
 for all $i\in\{1,\dots,q\}$. Hence, by \eqref{E:4.29+30}, we see that \eqref{E:4.30a} can be written as
  \begin{equation*}
   M\,\hw_0- L\mmatrix{(\vp^{[1]},\hw)_{N+1}\\ \vdots\\ (\vp^{[2q-2n]},\hw)_{N+1}}=0.
  \end{equation*}
 Therefore, the linear relation $T$ in \eqref{E:4.2} can be equivalently expressed as in \eqref{E:2.32}, which means that $T$ is 
 a self-adjoint extension of $\Tmin$. 
\end{proof}

The simplification of Theorem~\ref{T:s-e.ext} in the limit circle case is based on the following lemma.

\begin{lemma}\label{L:3.4}
 Let Hypothesis~\ref{H:SAC} be satisfied and $\vp^{[1]},\dots,\vp^{[\qp]}$ be arranged as in Lemma~\ref{L:3.2}. Assume that 
 there exists $\nu\in\Rbb$ such that system \Sla{\nu} has $r\coloneq\max\{\qp,\qm\}$ linearly independent square summable 
 solutions (suppressing the argument $\nu$) given by $\Th^{[1]},\dots,\Th^{[r]}$. Then these solutions can be arranged such 
 that $\rank\Ups_{p-2n}=p-2n$, where
  \begin{equation*}\label{E:3.24}
   \Ups\coloneq 
    \mmatrix{(\Th^{[1]},\Th^{[1]})_{N+1} & \hdots & (\Th^{[1]},\Th^{[r]})_{N+1}\\
               \vdots & \ddots & \vdots\\
              (\Th^{[r]},\Th^{[1]})_{N+1} & \hdots & (\Th^{[r]},\Th^{[r]})_{N+1}}\in\Cbb^{r\times r}.
  \end{equation*}
 Moreover, for any $\{\tz,\tf\}\in\Tmax$ the element $\hz$ can be uniquely expressed as
  \begin{equation*}\label{E:3.25}
   \hz_k=\hy_k+\sum_{i=1}^{2n} \al_{i}\,\hz_k^{[i]}+\sum_{j=1}^{p-2n} \be_j\,\Th_k^{[j]}, k\in\Izp,
  \end{equation*}
 where $\hy\in\dom\Tmin$, $\hz^{[1]},\dots,\hz^{[2n]}$ are given in Lemma~\ref{L:2.6}, and $\al_i,\be_j\in\Cbb$ for all 
 $i\in\{1,\dots,2n\}$ and $j\in\{1,\dots,p-2n\}$. 
\end{lemma}
\begin{proof}
 Since $\Th^{[1]},\dots,\Th^{[r]}\in\dom\Tmax$, by Lemma~\ref{L:3.2} there exist unique $\al_{i,j},\be_{i,l}\in\Cbb$ such that
  \begin{equation}\label{E:3.26}
   \Th_k^{[i]}=\hy^{[i]}_k+\sum_{j=1}^{2n} \al_{i,j}\,\hz_k^{[j]}+\sum_{l=1}^{p-2n} \be_{i,l}\,\vp_k^{[l]},\quad k\in\Izp,
  \end{equation}
 where $i\in\{1,\dots,r\}$. Then, the definition of $\hz^{[i]}$ and identity \eqref{E:2.39} yield
  \begin{equation}\label{E:3.27}
   \Ups=B\,\Om_{p-2n}\,B^*,
  \end{equation}
 where the matrix $B=[\ \overline{\be_{i,j}}\ ]\in\Cbb^{r\times(p-2n)}$. Hence, $\rank\Ups\leq p-2n$ by the first 
 inequality in \eqref{E:2.1}. On the other hand, by the Wronskian-type identity in \eqref{E:wronski.id} we have  
 $\Ups=\Th_0^*\,\Jc\,\Th_0$, where $\Th_k\coloneq(\Th^{[1]}_k,\dots,\Th_k^{[r]})$. Since the solutions  
 $\Th^{[1]}_k,\dots,\Th_k^{[r]}$ are linearly independent, we have $\rank\Th_k=r$ for all $k\in\Izp$, and hence
 $\rank\Ups\geq p-2n$ by the second inequality in \eqref{E:2.1}. Therefore $\rank\Ups=p-2n$, which implies that the solutions 
 $\Th_k\coloneq(\Th^{[1]}_k,\dots,\Th_k^{[r]})$ can be arranged such that $\rank\Ups_{p-2n}=p-2n$. In this case, the 
 invertibility of $B_{p-2n}$ follows from the equality $\Ups_{p-2n}=B_{p-2n}\,\Om_{p-2n}\,B^*_{p-2n}$, which is obtained 
 analogously to \eqref{E:3.27}. Since from \eqref{E:3.26} we have
  \begin{equation*}
   (\Th^{[1]}_k,\dots,\Th^{[p-2n]}_k)=(\hy^{[1]}_k,\dots,\hy^{[p-2n]}_k)+(\hz^{[1]}_k,\dots,\hz^{[2n]}_k)\,A^*_{2n,p-2n}
     +(\vp^{[1]}_k,\dots,\vp^{[p-2n]}_k)\,B^*_{q-2n},
  \end{equation*}
 where $A=[\ \overline{\al_{i,j}}\ ]\in\Cbb^{r\times2n}$, the invertibility of $B_{p-2n}$ means that  
 $\vp^{[1]}_k,\dots,\vp^{[p-2n]}_k$ can be uniquely expressed by using $\Th^{[1]}_k,\dots,\Th^{[p-2n]}_k$, 
 $\hy^{[1]}_k,\dots,\hy^{[p-2n]}_k$, and $\hz^{[1]}_k,\dots,\hz^{[2n]}_k$. Upon combining these expressions with \eqref{E:3.8}, 
 we obtain the second part of the statement.
\end{proof}

\section*{Acknowledgements}
This work was supported by the Program of ``Employment of Newly Graduated Doctors of Science for Scientific Excellence'' (grant 
number CZ.1.07/2.3.00/30.0009) co-financed from European Social Fund and the state budget of the Czech Republic. The first 
author would like to express his thanks to the Department of Mathematics and Statistics (Missouri University of Science 
and Technology) for hosting his visit. The authors are also indebted to the anonymous referee for detailed reading of 
the manuscript and constructive comments in her/his report which helped to improve the presentation of the results.



\end{document}